\documentclass[11pt]{amsart}

\usepackage[all]{xy}
\usepackage{color, pstricks}
\usepackage{mathrsfs}
\usepackage{amsmath,amsthm, amssymb, amsgen,amsxtra,amsfonts, amsbsy} 
\usepackage[all]{xy}
\usepackage{verbatim}

\newcommand{\mathds}[1]{{\mathbb #1}}

\newcommand{\diag}{\mathrm{diag}}
\newcommand{\nondiag}{\mathrm{nondiag}}

\newcommand{\pt}{\mathrm{pt}}
\DeclareMathOperator{\pr}{pr}
\DeclareMathOperator{\Gal}{Gal}

\DeclareMathOperator{\Sing}{Sing}
\DeclareMathOperator{\Ker}{Ker}

\begin{document}
%
%   D e f i n i t i o n s
%
%
\theoremstyle{definition}
\newtheorem{Definition}{Definition}[section]
\newtheorem*{Definitionx}{Definition}
\newtheorem{Convention}{Definition}[section]
\newtheorem{Construction}{Construction}[section]
\newtheorem{Example}[Definition]{Example}
\newtheorem{Examples}[Definition]{Examples}
\newtheorem{Remark}[Definition]{Remark}
\newtheorem*{Remarkx}{Remark}
\newtheorem{Remarks}[Definition]{Remarks}
\newtheorem{Caution}[Definition]{Caution}
\newtheorem{Conjecture}[Definition]{Conjecture}
\newtheorem*{Conjecturex}{Conjecture}
\newtheorem{Question}{Question}
\newtheorem{Questions}[Definition]{Questions}
\newtheorem*{Acknowledgements}{Acknowledgements}
\newtheorem*{AssumptionStar}{Assumption $(\star)$}
\newtheorem*{Organization}{Organization}
\newtheorem*{Disclaimer}{Disclaimer}
\theoremstyle{plain}
\newtheorem{Theorem}[Definition]{Theorem}
\newtheorem*{Theoremx}{Theorem}
\newtheorem{Proposition}[Definition]{Proposition}
\newtheorem*{Propositionx}{Proposition}
\newtheorem{Lemma}[Definition]{Lemma}
\newtheorem{Corollary}[Definition]{Corollary}
\newtheorem*{Corollaryx}{Corollary}
\newtheorem{Fact}[Definition]{Fact}
\newtheorem{Facts}[Definition]{Facts}
\newtheoremstyle{voiditstyle}{3pt}{3pt}{\itshape}{\parindent}%
{\bfseries}{.}{ }{\thmnote{#3}}%
\theoremstyle{voiditstyle}
\newtheorem*{VoidItalic}{}
\newtheoremstyle{voidromstyle}{3pt}{3pt}{\rm}{\parindent}%
{\bfseries}{.}{ }{\thmnote{#3}}%
\theoremstyle{voidromstyle}
\newtheorem*{VoidRoman}{}

% abgeschrieben aus The LaTeX Companion, 2nd edition,
% von Mittelback & Goossens
%
%\newcommand{\prf}{\par\noindent{\sc Proof.}\quad}
\newcommand{\prf}[1][Proof]{\par\noindent{\sc #1.}\quad}
\newcommand{\blowup}{\rule[-3mm]{0mm}{0mm}}
\newcommand{\cal}{\mathcal}
\newcommand{\Aff}{{\mathds{A}}}
\newcommand{\BB}{{\mathds{B}}}
\newcommand{\CC}{{\mathds{C}}}
\newcommand{\EE}{{\mathds{E}}}
\newcommand{\FF}{{\mathds{F}}}
\newcommand{\GG}{{\mathds{G}}}
\newcommand{\HH}{{\mathds{H}}}
\newcommand{\NN}{{\mathds{N}}}
\newcommand{\ZZ}{{\mathds{Z}}}
\newcommand{\PP}{{\mathds{P}}}
\newcommand{\QQ}{{\mathds{Q}}}
\newcommand{\RR}{{\mathds{R}}}
\newcommand{\Sphere}{{\mathds{S}}}
\newcommand{\Liea}{{\mathfrak a}}
\newcommand{\Lieb}{{\mathfrak b}}
\newcommand{\Lieg}{{\mathfrak g}}
\newcommand{\Liem}{{\mathfrak m}}
\newcommand{\ideala}{{\mathfrak a}}
\newcommand{\idealb}{{\mathfrak b}}
\newcommand{\idealg}{{\mathfrak g}}
\newcommand{\idealm}{{\mathfrak m}}
\newcommand{\idealp}{{\mathfrak p}}
\newcommand{\idealq}{{\mathfrak q}}
\newcommand{\idealI}{{\cal I}}
\newcommand{\lin}{\sim}
\newcommand{\num}{\equiv}
\newcommand{\dual}{\ast}
\newcommand{\iso}{\cong}
\newcommand{\homeo}{\approx}
\newcommand{\mm}{{\mathfrak m}}
\newcommand{\pp}{{\mathfrak p}}
\newcommand{\qq}{{\mathfrak q}}
\newcommand{\rr}{{\mathfrak r}}
\newcommand{\pP}{{\mathfrak P}}
\newcommand{\qQ}{{\mathfrak Q}}
\newcommand{\rR}{{\mathfrak R}}
%
%  evtl. auch \"uber \mathbb oder \Bbb
%
\newcommand{\OO}{{\cal O}}
\newcommand{\numero}{{n$^{\rm o}\:$}}
\newcommand{\mf}[1]{\mathfrak{#1}}
\newcommand{\mc}[1]{\mathcal{#1}}
\newcommand{\into}{{\hookrightarrow}}
\newcommand{\onto}{{\twoheadrightarrow}}
\newcommand{\Spec}{{\rm Spec}\:}
\newcommand{\BigSpec}{{\rm\bf Spec}\:}
\newcommand{\Spf}{{\rm Spf}\:}
\newcommand{\Proj}{{\rm Proj}\:}
\newcommand{\Pic}{{\rm Pic }}
\newcommand{\Br}{{\rm Br}}
\newcommand{\NS}{{\rm NS}}
\newcommand{\Sym}{{\mathfrak S}}
\newcommand{\Aut}{{\rm Aut}}
\newcommand{\Autp}{{\rm Aut}^p}
\newcommand{\Hom}{{\rm Hom}}
\newcommand{\Ext}{{\rm Ext}}
\newcommand{\ord}{{\rm ord}}
\newcommand{\coker}{{\rm coker}\,}
\newcommand{\id}{{\rm id}}
\newcommand{\divisor}{{\rm div}}
\newcommand{\Def}{{\rm Def}}
\newcommand{\piet}{{\pi_1^{\rm \acute{e}t}}}
\newcommand{\Het}[1]{{H_{\rm \acute{e}t}^{{#1}}}}
\newcommand{\Pet}[1]{{P_{\rm \acute{e}t}^{{#1}}}}
\newcommand{\Hcris}[1]{{H_{\rm cris}^{{#1}}}}
\newcommand{\HdR}[1]{{H_{\rm dR}^{{#1}}}}
\newcommand{\hdR}[1]{{h_{\rm dR}^{{#1}}}}
\newcommand{\defin}[1]{{\bf #1}}

\newcommand{\X}{{\mathcal X}}
\newcommand{\Z}{{\mathcal Z}}

\title{Good Reduction of K3 Surfaces}

\author{Christian Liedtke}
\address{TU M\"unchen, Zentrum Mathematik - M11, Boltzmannstr. 3, D-85748 Garching bei M\"unchen, Germany}
\email{liedtke@ma.tum.de}
\author{Yuya Matsumoto}
\address{Graduate School of Mathematics, Nagoya University, Furocho, Chikusaku, Nagoya, 464-8602, Japan}
\email{matsumoto.yuya@math.nagoya-u.ac.jp}

\date{March 6, 2017}
% K3, varieties over local fields, Galois representations, MMP and Mori theory 

\begin{abstract}
  Let $K$ be the field of fractions of a local Henselian 
  discrete valuation ring
  $\OO_K$ of characteristic zero with perfect residue field $k$.
  Assuming potential semi-stable reduction, we show 
  that an unramified Galois action on the second $\ell$-adic cohomology group
  of a K3 surface over $K$ implies that the 
  surface has good reduction after a finite and {\em unramified} extension.
  We give examples where this unramified extension
  is really needed.
  Moreover, we give applications to good reduction after tame extensions
  and Kuga--Satake Abelian varieties.
  On our way, we settle existence and termination of
  certain flops in mixed characteristic, and study group actions
  and their quotients on models of varieties.
\end{abstract}

\maketitle

\section{Introduction}

Let $\OO_K$ be a local Henselian discrete valuation ring
of characteristic zero with
field of fractions $K$ and perfect residue field $k$, whose characteristic
is $p\geq0$.
For example, $\OO_K$ could be $\CC[[t]]$ or the ring of integers
in a $p$-adic field.
Given a variety $X$ that is smooth and proper over $K$, one can ask whether
$X$ has {\em good reduction}, that is, whether there exists an algebraic space
$$
  {\cal X}\,\to\,\Spec\OO_K
$$
with generic fiber $X$ that is smooth and proper over $\OO_K$.

\subsection{Good reduction and Galois representations}
Let $\ell$ be a prime different from $p$, let $G_K:={\rm Gal}(\overline{K}/K)$ 
be the absolute Galois group of $K$, and let $I_K$ be its inertia subgroup.
Then, the natural $\ell$-adic Galois representation
$$
 \begin{array}{ccccc}
   \rho_{m,\ell} &:& G_K &\to& {\rm Aut} \, \left( \, \Het{m}(X_{\overline{K}},\QQ_\ell) \,\right)
  \end{array}
$$
is called {\em unramified} if it satisfies $\rho_{m,\ell}(I_K)=\{{\rm id}\}$.
A necessary condition for $X$ to have good reduction is that for all $m\geq1$
and all primes $\ell\neq p$, the representation $\rho_{m,\ell}$ is 
unramified.

\subsection{Curves and Abelian varieties}
By a famous theorem of Serre and Tate \cite{serre tate}, 
which generalizes results of N\'eron, Ogg, and
Shafarevich for elliptic curves to Abelian varieties, 
the $G_K$-representation $\rho_{1,\ell}$ detects the 
reduction type of Abelian varieties.

\begin{Theorem}[Serre--Tate]
   An Abelian variety $X$ over $K$ 
   has good reduction if and only if the $G_K$-representation on
   $\Het{1}(X_{\overline{K}},\QQ_\ell)$ is unramified.
\end{Theorem}

On the other hand, it is not too difficult to give counterexamples to such a result
for curves of genus $\geq2$.
Nevertheless, Oda \cite{Oda} showed that good reduction can be detected by 
the outer $G_K$-representation on the \'etale fundamental group.
We refer the interested reader to Section \ref{subsec: counterexamples}
for references, examples, and details.

\subsection{Kulikov--Nakkajima--Persson--Pinkham models}
Before coming to the results of this article, we have to make one crucial
assumption.

\begin{AssumptionStar}
  A K3 surface $X$ over $K$ satisfies ($\star$) if
  there exists a finite field extension $L/K$ such that
  $X_L$ admits a model ${\cal X}\to\Spec\OO_L$
  that is a regular algebraic space
  with trivial canonical sheaf $\omega_{{\cal X}/\OO_L}$, and 
  whose geometric special fiber is a normal crossing divisor.
\end{AssumptionStar}

In equal characteristic zero, ($\star$) always holds,
and the special fibers of the corresponding models 
have been classified by Kulikov \cite{Kulikov}, Persson \cite{Persson}, and Persson--Pinkham \cite{PP}.
In mixed characteristic, the corresponding classification 
(assuming the existence of such models) is due to Nakkajima \cite{Nakkajima}.
If the expected results on resolution of singularities and toroidalization of morphisms
were known to hold in mixed characteristic, then $(\star)$ would follow from
Kawamata's semi-stable minimal model program (MMP) 
in mixed characteristic \cite{Kawamata} and Artin's results \cite{Artin Brieskorn} 
on simultaneous resolutions of families of surface singularities.
We refer to Proposition \ref{prop: star is satisfied} for details.
Using work of Maulik \cite{Maulik} and some strengthenings 
due to the second named author \cite{matsumoto2}, we have at least the following.

\begin{Theorem}
  \label{thm: star holds}
  Let $X$ be a K3 surface over $K$ and assume that $p=0$ or that 
  $X$ admits an ample invertible sheaf ${\cal L}$ with $p>{\cal L}^2+4$.
  Then, $X$ satisfies $(\star)$.
\end{Theorem}

\subsection{K3 surfaces}

In this article, we establish a N\'eron--Ogg--Shafarevich--Serre--Tate 
type result for K3 surfaces.
Important steps were already taken by the second named author 
in \cite{matsumoto2}.
Over the complex numbers, similar results are classically known, see, 
for example, \cite[Chapter 5]{Kulikov survey}.

Before coming to the main result of this article, we 
define a \emph{K3 surface with at worst RDP singularities} 
to be a proper surface over a field, which, after base change to an algebraically closed field,
has at worst rational double point singularities, and whose minimal resolution 
of singularities is a K3 surface.

\begin{Theorem}
 \label{thm: main in introduction}
  Let $X$ be a K3 surface over $K$ that satisfies $(\star)$.
  If the $G_K$-representation on $\Het{2}(X_{\overline{K}},\QQ_\ell)$ is
  unramified for some $\ell\neq p$, then
  \begin{enumerate}
    \item there exists a model of $X$ that is a projective scheme over $\OO_K$, 
       whose special fiber is a K3 surface with at worst RDP singularities.
    \item Moreover, there exists an integer $N$, independent of $X$ and $K$, and a
       finite unramified extension $L/K$ of degree $\leq N$, such that
       $X_L$ has good reduction over $L$.
  \end{enumerate}
\end{Theorem}

In \cite[Theorem 35]{Hassett Tschinkel}, a similar result is obtained for K3 surfaces
over $\CC((t))$, but their proof uses methods different from ours.
As in the case of Abelian varieties in \cite{serre tate}, we obtain the
following independence of $\ell$.

\begin{Corollary}
  Let $X$ be a K3 surface over $K$ that satisfies $(\star)$.
  Then, the $G_K$-representation on $\Het{2}(X_{\overline{K}},\QQ_\ell)$
  is unramified for one $\ell\neq p$ if and only if it is unramified for all $\ell\neq p$.
\end{Corollary} 

In \cite{serre tate}, Serre and Tate showed that if an Abelian variety of dimension $g$
over $K$ with $p>2g+1$ has potential good reduction, then
good reduction can be achieved after a {\em tame} extension.
Here, we establish the following analog for K3 surfaces.

\begin{Corollary}
  Let $X$ be a K3 surface over $K$ with $p\geq23$ and potential good reduction.
  Then, $X$ has good reduction after a tame extension of $K$.
\end{Corollary}

It is important to note that in part (2) of Theorem \ref{thm: main in introduction}, 
we cannot avoid field extensions in general.
More precisely, we construct the following explicit examples.

\begin{Theorem}
  For every prime $p\geq5$, there exists a K3 surface $X=X(p)$ over $\QQ_p$, such that
  \begin{enumerate}
    \item the $G_{\QQ_p}$-representation on $\Het{2}(X_{\overline{\QQ}_p},\QQ_\ell)$ is unramified for all $\ell\neq p$, 
    \item $X$ has good reduction over the unramified extension $\QQ_{p^2}$, but
    \item $X$ does not have good reduction over $\QQ_p$.
  \end{enumerate}
\end{Theorem}

\subsection{Kuga--Satake Abelian varieties}
Let us recall that Kuga and Satake \cite{Kuga Satake} associated 
to a polarized K3 surface $(X,{\cal L})$ over $\CC$ a polarized Abelian
variety ${\rm KS}(X,{\cal L})$ of dimension $2^{19}$ over $\CC$.
Moreover, if $(X,{\cal L})$ is defined over an arbitrary field $k$, then
Rizov \cite{Rizov} and Madapusi Pera \cite{madapusi pera},
building on work of Deligne \cite{Deligne K3} and Andr\'e \cite{Andre},
established the existence of ${\rm KS}(X,{\cal L})$ over 
some finite extension of $k$.
As an application of Theorem \ref{thm: main in introduction},
we can compare the reduction behavior of a polarized
K3 surface to that of its associated Kuga--Satake Abelian variety.

\begin{Theorem}
  Assume $p \neq 2$.
  Let $(X,{\cal L})$ be a polarized K3 surface over $K$.
   \begin{enumerate}
     \item If $X$ has good reduction, then ${\rm KS}(X,{\cal L})$
       can be defined over an unramified extension $L/K$, and it has good reduction over $L$.
     \item Assume that $X$ satisfies $(\star)$. Let
       $L/K$ be a field extension such that 
       both ${\rm KS}(X,{\cal L})$ and the Kuga--Satake correspondence can be 
       defined over $L$.
       If ${\rm KS}(X,{\cal L})$ has good reduction over $L$, then $X$ has good reduction
       over an unramified extension of $L$.
  \end{enumerate}
\end{Theorem}

\subsection{Organization}
This article is organized as follows:

In Section \ref{sec: generalities}, we recall a couple of general facts on 
models and unramified Galois representations on $\ell$-adic cohomology. 
We also recall the classical Serre--Tate theorem for Abelian varieties and give explicit
examples of curves of genus $\geq2$, where the Galois representation does not detect
bad reduction.

In Section \ref{sec: K3 and Enriques}, we review potential semi-stable reduction of K3 surfaces,
Kawamata's semi-stable MMP,
the Kulikov--Nakkajima--Pinkham--Persson classification list, and the the second named 
author's results on potential good reduction of K3 surfaces.
% in terms of unramified Galois representations.
We also briefly discuss potential good and semi-stable reduction 
of Enriques surfaces.

In Section \ref{sec: flops}, we establish existence and termination of certain 
flops, which we need later on to equip our models with suitable invertible sheaves.
Moreover, we show that any two smooth models of a K3 surface $X$ over $K$ are related
by a finite sequence of flopping contractions and their inverses.

Section \ref{sec: group actions on models} is the technical heart of this article:
given a K3 surface $X$ over $K$, a finite Galois extension $L/K$ with group $G$, 
and a model of $X_L$ over $\OO_L$, 
we study extensions of the $G$-action $X_L$ to this model.
Then, we study quotients of such models by $G$-actions, where the most difficult
case arises when $p$ divides the order of $G$ (wild action).

In Section \ref{sec: main}, we establish the main results of this article:
a N\'eron--Ogg--Shafarevich type theorem for K3 surfaces, good reduction
over tame extensions, as well as the connection to Kuga--Satake Abelian
varieties.

Finally, in Section \ref{sec: counterexample}, we give explicit examples of K3 surfaces
over $\QQ_p$ with unramified Galois representations on their
$\ell$-adic cohomology groups that do not have good reduction over $\QQ_p$.
% have good reduction over an  unramified extension of $\QQ_p$, but not over $\QQ_p$ itself.

\begin{Acknowledgements}
 It is a pleasure for us to thank
 Fran\c{c}ois Charles, Christopher Hacon, Annabelle Hartmann, Brendan Hassett,
 Tetsushi Ito, Keerthi Madapusi Pera, Davesh Maulik, as well as Yoichi Mieda
 for discussions and comments.
 We also thank the referee for careful proofreading and many comments, 
 including an idea to simplify the arguments in Section \ref{sec: termination of flops}.
 The second named author thanks the department of mathematics at the TU M\"unchen 
 for kind hospitality while visiting the first named author there.
 The first named author is supported by the ERC Consolidator Grant
 681838 ``K3CRYSTAL'' and the
 second named author is supported by JSPS KAKENHI Grant Number 26247002.
\end{Acknowledgements}

\section*{Notations and Conventions}

Throughout the whole article, we fix the following notations
$$
 \begin{array}{ll} 
   \OO_K & \mbox{ a local Henselian DVR of characteristic zero } \\
   K & \mbox{ its field of fractions} \\
   k &\mbox{ the residue field, which we assume to be perfect}\\
   p\geq0 &\mbox{ the characteristic of $k$} \\
      \ell & \mbox{ a prime different from $p$} \\
   G_K,\, G_k & \mbox{ the absolute Galois groups }{\rm Gal}(\overline{K}/K), {\rm Gal}(\overline{k}/k) 
 \end{array}
$$
%By a variety over a field $F$, we mean a geometrically integral scheme of finite type over $F$.
If $L/K$ is a field extension, and $X$ is a scheme over $K$, we abbreviate
the base-change $X\times_{\Spec K}\Spec L$ by $X_L$.

\section{Generalities}
\label{sec: generalities}

In this section, we recall a couple of general facts on models of varieties,
unramified Galois representations on $\ell$-adic cohomology groups, 
and N\'eron--Ogg--Shafarevich type theorems.

\subsection{Models} 
We start with the definition of various types of models.

\begin{Definition}
   \label{def:model}
    Let $X$ be a smooth and proper variety over $K$.
    \begin{enumerate}
      \item  A {\em model of $X$ over $\OO_K$} is an algebraic space that is flat and 
         proper over $\Spec\OO_K$ and whose generic fiber is isomorphic to $X$.
      \item We say that $X$ has {\em good reduction} if there exists a model of $X$ that
         is smooth over $\OO_K$.
      \item We say that $X$ has {\em semi-stable reduction} if there exists a regular
        model of $X$, whose geometric special fiber is a reduced normal crossing divisor
        with smooth components.
        (Sometimes, this notion is also called \emph{strictly} semi-stable reduction.)
      \item We say that $X$ has {\em potential} good (resp. semi-stable) reduction
        if there exists a finite field extension $L/K$ such that $X_L$ has good 
        (resp. semi-stable) reduction.
   \end{enumerate}
\end{Definition}

\begin{Remark}
  \label{rem: curves an Abelian varieties}
  Models of curves and Abelian varieties can be treated entirely within the category of
  schemes, see, for example, \cite[Chapter 10]{Liu} and \cite{BLR}.
  However, if $X$ is a K3 surface over $K$ with good reduction, then it may
  not be possible to find a smooth model in the category of schemes, and
  we refer to \cite[Section 5.2]{matsumoto2} for explicit examples.
  In particular, we are forced to work with algebraic spaces when studying models
  of K3 surfaces.
\end{Remark}  

\subsection{Inertia and monodromy}
\label{subsec: inertia}
%We recall, for example from \cite[Chapter I.2]{serre galois} %\cite[Section I.7]{serre corps},
The $G_K$-action on $\overline{K}$ induces an action on $\OO_{\overline{K}}$
and by reduction, an action on $\overline{k}$.
This gives rise to a continuous and surjective homomorphism $G_K\to G_k$ 
of profinite groups. 
Thus, we obtain a short exact sequence
$$
1\,\to\,I_K\,\to\,G_K\,\to\,G_k\,\to\,1,
$$
whose kernel $I_K$ is called the {\em inertia group}.
In fact, $I_K$ is the absolute Galois group 
of the maximal unramified extension of $K$.
If $p\neq0$, then the {\em wild inertia group} $P_K$ is the normal subgroup 
of $G_K$ that is the absolute Galois group of the maximal 
tame extension of $K$.
We note that $P_K$ is the unique $p$-Sylow subgroup of $I_K$.

\begin{Definition}
 \label{def: unramified action}
  Let $X$ be a smooth and proper variety over $K$.
  Then, the $G_K$-representation on $\Het{m}(X_{\overline{K}},\QQ_\ell)$ is 
  called {\em unramified} if $I_K$ acts trivially.
  It is called {\em tame} if $P_K$ acts trivially.
\end{Definition}

For an Abelian variety $X$, it follows from results of Serre and Tate \cite{serre tate} 
that the $G_K$-representation on $\Het{m}(X_{\overline{K}},\QQ_\ell)$ is unramified
for {\em one} $\ell\neq p$, if and only if it is so for {\em all} $\ell\neq p$.
In Corollary \ref{cor: independence of l}, we will show a similar result for K3 surfaces.
In general, it is not known whether being unramified depends on the choice of $\ell$,
but it is expected not to.

A relation between good reduction and unramified
Galois representations on $\ell$-adic cohomology groups is given
by the following well-known result, which follows from the proper smooth 
base change theorem.
For schemes, it is stated in \cite[Th\'eor\`eme XII.5.1]{SGA4}, 
and in case the model is an algebraic space,
we refer to \cite[Theorem 0.1.1]{Liu Zheng} 
or \cite[Chapitre VII]{Artin representabilite}.

\begin{Theorem}
 \label{thm:obvious}
  If $X$ has good reduction, then the $G_K$-representation
  on $\Het{m}(X_{\overline{K}},\QQ_\ell)$ is unramified for all $m$ and for all 
  $\ell\neq p$.
\end{Theorem}

In view of this theorem, it is natural to ask for the converse direction.
Whenever such a converse holds for some class of varieties over $K$, 
we obtain a purely representation-theoretic
criterion to determine whether such a variety admits a model 
over $\OO_K$ with good reduction.

\subsection{Abelian varieties}
A classical converse to Theorem \ref{thm:obvious} is the N\'eron--Ogg--Shafarevich 
criterion for elliptic curves.
Later, Serre and Tate generalized it to Abelian varieties of arbitrary dimension.

\begin{Theorem}[Serre--Tate \cite{serre tate}]
 \label{thm:serre tate}
 An Abelian variety $A$ over $K$ has good reduction if and only if the
 $G_K$-representation on $\Het{1}(A_{\overline{K}},\QQ_\ell)$ is unramified
 for one (resp. all) $\ell\neq p$.
\end{Theorem}

\subsection{Higher genus curves, part 1}
\label{subsec: counterexamples}
Now, the converse to Theorem \ref{thm:obvious} already fails for curves of higher genus.
Let $X$ be a smooth and proper curve of genus $g\geq2$ over $K$.
Let ${\rm Jac}(X)$ be its Jacobian, which is an Abelian variety of dimension $g$ over $K$.
Then, the exact sequence of \'etale sheaves on $X$
$$
 1 \,\to\, \mu_n\, \to\, \GG_m \,\stackrel{\times n}{\longrightarrow}\, \GG_m \,\to\, 1
$$
%gives rise to a $G_K$-equivariant isomorphism between
%$\Het{1}({\rm Jac}(X)_{\overline{K}},\QQ_\ell)$ and $\Het{1}(X_{\overline{K}},\QQ_\ell)$.
gives rise to $G_K$-equivariant isomorphisms
$\Het{1}(X_{\overline{K}},\mu_n) \cong \Pic(X_{\overline{K}})[n] \cong \Het{1}({\rm Jac}(X)_{\overline{K}},\mu_n)$, 
from which we obtain
$\Het{1}(X_{\overline{K}},\QQ_\ell) \cong \Het{1}({\rm Jac}(X)_{\overline{K}},\QQ_\ell)$
by passing to the limit.
Moreover, if $X$ has a $K$-rational point, then there is a natural embedding
$$
   j\,:\, X \,\to\, {\rm Jac}(X),
$$
and the above isomorphism coincides with $j^*$ (which is independent of the choice of a rational point).
By the Serre--Tate theorem (Theorem \ref{thm:serre tate}), 
an unramified $G_K$-representation  on $\Het{1}(X_{\overline{K}},\QQ_\ell)$ 
is equivalent to good reduction of ${\rm Jac}(X)$.
The following lemma gives a criterion that ensures the latter.

\begin{Lemma}
  \label{lem: Jacobian good reduction}
  Let $X$ be a smooth and proper curve over $K$ that admits a semi-stable 
  scheme model ${\cal X}\to \Spec\OO_K$ such that the dual graph associated 
  to the components of its special fiber ${\cal X}_0$ is a tree.
  Then, ${\rm Jac}(X)$ has good reduction and the $G_K$-representation on 
  $\Het{1}(X_{\overline{K}},\QQ_\ell)$ is unramified.
\end{Lemma}

\prf
By \cite[Section 9.2, Example 8]{BLR},  $\Pic^0_{{\cal X}_0/k}$ is an Abelian
variety, which implies that ${\rm Jac}(X)$ has good reduction, and thus,
the $G_K$-representation on $\Het{1}(X_{\overline{K}},\QQ_\ell)$
is unramified.
\qed\medskip

Using this lemma, it is easy to produce counterexamples
to N\'eron--Ogg--Shafarevich type results for curves of higher genus.

\begin{Proposition}
  If $p\neq2$, then there exists for infinitely many $g\geq2$ a smooth
  and proper curve $X$ of genus $g$ over $K$ such that
  \begin{enumerate}
   \item the $G_K$-representation on $\Het{m}(X_{\overline{K}},\QQ_\ell)$
    is unramified for all $m$ and all $\ell\neq p$, and 
   \item $X$ does not have good reduction over $K$ nor over any finite extension.
  \end{enumerate}
\end{Proposition}

\prf
We give examples for $g \not\equiv 1\mod p$.
Let $X$ be a hyperelliptic curve of genus $g$ over $K$ that is one of the examples of
\cite[Example 10.1.30]{Liu} with the extra assumptions of \cite[Example 10.3.46]{Liu} 
(here, we need the assumption $g \not\equiv 1 \mod p$).
% $X$ should furthermore satisfy the condition $\disc(X^{2g-1} + c_2 X + d_2) \in \OO_K^*$,
% and we leave it to the interested reader to fill the details.
Then, $X$ has stable reduction over $K$, as well as over every finite extension field $L/K$.
In this example, the special fiber of the stable model is the union of 
a curve of genus $1$ and a curve of  genus $(g-1)$ meeting transversally in one point.
In particular, neither $X$ nor any base-change $X_L$ have good reduction, but since the assumptions of
Lemma \ref{lem: Jacobian good reduction} are fulfilled, the $G_K$-representation on 
$\Het{m}(X_{\overline{K}},\QQ_\ell)$ is unramified for all $m$ and all $\ell\neq p$.
\qed\medskip

We stress that these results are well-known to the experts, but since we were not able to find explicit 
references and explicit examples, we decided to include them here.

\subsection{Higher genus curves, part 2}
If $X$ is a smooth and proper curve of genus $\geq2$ over $K$, then one can also study the outer 
$G_K$-representation on its \'etale fundamental group, 
which turns out to detect good reduction.
More precisely, there exists a short exact sequence of \'etale fundamental groups
$$
  1\,\to\,\piet\left( X_{\overline{K}} \right)\,\to\,\piet\left(X\right)\,\to\,G_K \,\to\,1\,.
$$
For every prime $\ell$, 
this exact sequence gives rise to a well-defined homomorphism from $G_K$ to the 
the outer automorphism group of the pro-$\ell$-completion $\piet(X_{\overline{K}})_\ell$
of the geometric \'etale fundamental group
$$
  \begin{array}{ccccc}
     \rho_\ell &:& G_K &\longrightarrow& {\rm Out}\left( \piet(X_{\overline{K}})_\ell  \right) \,. 
   \end{array}
$$
In analogy to Definition \ref{def: unramified action},
we will say that this representation is {\em unramified} if $\rho_\ell(I_K)=\{1\}$.
We note that the $G_K$-representation on $\Het{1}(X_{\overline{K}},\QQ_\ell)$ arises from
the residual action of $\rho_\ell$ on the Abelianization of $\piet(X_{\overline{K}})_\ell$. 
After these preparations, we have the following N\'eron--Ogg--Shafarevich type 
theorem for curves of higher genus, which is in terms of  fundamental groups rather
than cohomology groups.

\begin{Theorem}[{Oda \cite[Theorem 3.2]{Oda}}]
  Let $X$ be a smooth and proper curve of genus $\geq2$ over $K$.
  Then, $X$ has good reduction if and only if the outer Galois action $\rho_\ell$ 
  is unramified for one (resp. all) $\ell\neq p$.
\end{Theorem}

\section{K3 surfaces and their models}
\label{sec: K3 and Enriques}

In this section, we first introduce the crucial Assumption $(\star)$, which ensures the
existence of suitable models for K3 surfaces.
These models have been studied by Kulikov, Nakkajima, Persson, and Pinkham.
Following ideas of Maulik, 
we show how $(\star)$ would follow from a combination of 
potential semi-stable reduction (which is not known in mixed characteristic, 
but expected) and the semi-stable minimal model program (MMP) in mixed characteristic.
Then, we give some conditions under which $(\star)$ does hold.
After that, we shortly review the second named author's results on potential
good reduction of K3 surfaces.
Finally, we show by example that these results do not carry over to Enriques surfaces.
Most of the results of this section are probably known to the experts.

\subsection{Kulikov--Nakkajima--Persson--Pinkham models}
\label{subsec: potential semi-stable reduction}
We first introduce the crucial assumption that we shall make from now on.

\begin{AssumptionStar}
  A K3 surface $X$ over $K$ satisfies ($\star$) if
  there exists a finite field extension $L/K$ such that
  $X_L$ admits a semi-stable model ${\cal X}\to\Spec\OO_L$
  (in the sense of Definition \ref{def:model})
  such that $\omega_{{\cal X}/\OO_L}$ is trivial.
\end{AssumptionStar}

Here, we equip $\cal X$ with its standard log structure 
${\cal X}^{\mathrm{log}}$ and define the relative canonical sheaf
$\omega_{{\cal X}/\OO_L}$ to be 
$\bigwedge^2 \Omega^1_{{\cal X}^{\mathrm{log}}/\OO_L^{\mathrm{log}}}$ 
using log differentials.
Since ${\cal X}^{\mathrm{log}}$ is log smooth over $\OO_L$, 
the sheaf $\omega_{{\cal X}/\OO_L}$ is invertible,
see also the discussion in \cite[Section 3]{matsumoto2}.
%In the more general case appearing in the proof of the next proposition 
%we define $\omega_{{\cal X}/\OO_L} = H^{-2}()$,
%which is a reflexive sheaf of rank $1$. 
%In the semistable case the two definitions agree.
%****

The main reason why $(\star)$ is not known to hold is that potential
semi-stable reduction is not known:
using resolution of singularities
in mixed characteristic (recently announced by Cossart--Piltant \cite{Cossart--Piltant})
and embedded resolution of singularities (Cossart--Jannsen--Saito \cite{Cossart--Jannsen--Saito}),
we obtain a model $\cal X$, whose special fiber ${\cal X}_0$ 
has simple normal crossing support, but whose components may have multiplicities.
At the moment, it is not clear how to get rid of these multiplicities after base change, 
unless all of them are prime to $p$. 
In case the residue characteristic is zero, these results are classically known to hold, 
see the discussion in \cite[Section 7.2]{kollar mori} for details.

The following result, which is inspired by Maulik's approach and ideas 
from \cite[Section 4]{Maulik},
shows that $(\star)$ essentially holds once we assume potential semi-stable reduction.
More precisely, we have the following.

\begin{Proposition}
  \label{prop: star is satisfied}
   Assume $p \neq 2,3$. Let $X$ be a K3 surface over $K$ and
   assume that there exists 
   \begin{enumerate}
     \item a finite field extension $L'/K$, and
     \item a smooth surface $Y$ over $L'$ that is birationally equivalent to $X_{L'}$, and
     \item a scheme model ${\cal Y}\to\Spec\OO_{L'}$ of $Y$ with semi-stable reduction.
  \end{enumerate}
  Then, $X$ satisfies $(\star)$.
\end{Proposition}

\prf
Let ${\cal Y}\to\Spec\OO_{L'}$ be as in the statement.
Since $p\neq2,3$, Kawamata's semi-stable MMP \cite{Kawamata}
(see also \cite[Section 7.1]{kollar mori} for $p=0$)
produces a scheme ${\cal Z}\to\Spec\OO_{L'}$ with nef
relative canonical divisor $K_{{\cal Z}/\OO_{L'}}$ 
that is a model of a smooth proper surface birationally equivalent 
to $X_{L'}$, and such that ${\cal Z}$ 
is regular outside a finite set $\Sigma$ of terminal singularities.
We refer to loc.\ cit.\ for details and the definition of $K_{{\cal Z}/\OO_{L'}}$,
which is a Weil divisor.
We also note that it coincides with the Weil divisor class
associated to the relative canonical divisor $\omega_{{\cal Z}/\OO_{L'}}$, 
see, for example, \cite[Section 3]{matsumoto2}.

Since $X_{L'}$ is a minimal surface and $K_{{\cal Z}/\OO_{L'}}$ is nef, 
the generic fiber of ${\cal Z}$ is actually isomorphic to $X_{L'}$,
and it follows that $K_{{\cal Z}/\OO_{L'}}$ is trivial.
Outside $\Sigma$, this model is already a semi-stable model.
From the classification of terminal singularities in \cite[Theorem 4.4]{Kawamata}
and the fact that $K_{{\cal Z}/\OO_{L'}}$ is Cartier at points of $\Sigma$ (since it is trivial), 
it follows that the geometric special fiber $({\cal Z}_0)_{\overline{k}}$ 
is irreducible around points of $\Sigma$, and that it acquires RDP singularities in these points.
Thus, after some finite field extension $L/L'$, there exists a simultaneous
resolution ${\cal X}\to\Spec\OO_L$ of these singularities by \cite[Theorem 2]{Artin Brieskorn}.
This ${\cal X}$ may exist only as an algebraic space, and it satisfies $(\star)$.
\qed\medskip

As already mentioned above, the assumptions are fulfilled if $p=0$,
see \cite[Chapter 2]{Kempf} or the discussion in \cite[Section 7.2]{kollar mori}.
If $p\neq0$, then they are fulfilled for K3 surfaces that admit a very ample invertible
sheaf ${\cal L}$ with $p>{\cal L}^2+4$ by a result of Maulik \cite[Section 4]{Maulik}.
With some extra work, the condition ``very ample'' can be weakened to ``ample'' 
(see \cite[argument following Lemma 3.1]{matsumoto2})
and Theorem \ref{thm: star holds} follows.
Thus, we have the following result.

\begin{Theorem} (=Theorem \ref{thm: star holds})
  Let $X$ be a K3 surface over $K$ and assume that $p=0$ or that 
  $X$ admits an ample invertible sheaf ${\cal L}$ with $p>{\cal L}^2+4$.
  Then, $X$ satisfies $(\star)$.
\end{Theorem}

Over $\CC$, Kulikov \cite{Kulikov}, Persson \cite{Persson}, and Pinkham--Persson \cite{PP}
classified the special fibers of the models asserted by $(\star)$.
We refer to \cite[Section 1]{Morrison} and \cite[Chapter 5]{Kulikov survey}
for overview, and to Nakkajima's extension \cite{Nakkajima} 
of these results to mixed characteristic.

\subsection{Potential good reduction of K3 surfaces}
Now, if $X$ is a K3 surface over $K$ that satisfies $(\star)$, then
there exists a finite field extension $L/K$ and a model ${\cal X}\to\Spec\OO_L$
of $X_L$ as asserted by $(\star)$.
If the $G_K$-representation on $\Het{2}(X_{\overline{K}},\QQ_\ell)$ is unramified,
then the weight filtration on $\Het{2}(X_{\overline{K}},\QQ_\ell)$
that arises from the Steenbrink--Rapoport--Zink spectral sequence
(see \cite{Steenbrink}, \cite[Satz 2.10]{rapo}, and \cite[Proposition 1.9]{Nakayama} for details) 
is trivial.
Together with a result of Persson \cite[Proposition 3.3.6]{Persson},
this implies that the special fiber of ${\cal X}$ is smooth, that is,
$X_L$ has good reduction.
Thus, we obtain the following result of the second named author
and we refer to \cite{matsumoto2} for details and a detailed proof.

\begin{Theorem}[Matsumoto]
  \label{thm: matsumoto}
  Let $X$ be a K3 surface over $K$ that satisfies $(\star)$.
  If the $G_K$-representation on $\Het{2}(X_{\overline{K}},\QQ_\ell)$ 
  is unramified for one $\ell\neq p$, then $X$ has potential good reduction. 
\end{Theorem}

\subsection{Enriques surfaces}
The previous theorem does not generalize to other classes of
surfaces with numerically trivial canonical sheaves.
For example, the $G_K$-representation on $\ell$-adic cohomology
of an Enriques surface can neither exclude nor confirm any
type in the Kulikov--Nakkajima--Persson--Pinkham list for these surfaces.
More precisely, we have the following.

\begin{Lemma}
   \label{lem: enriques}
  Let $Y$ be an Enriques surface over $K$. 
  Then, there exists a finite extension $L/K$ such 
  that the $G_L$-representation on $\Het{m}(Y_{\overline{K}},\QQ_\ell)$ is 
  unramified for all $m$ and all $\ell\neq p$.
\end{Lemma} 

\prf
We only have to show something for $m=2$.
But then, the first Chern class induces a $G_K$-equivariant isomorphism
$$
  \NS(Y_{\overline{K}})\,\otimes_\ZZ\QQ_\ell\,\stackrel{c_1}{\longrightarrow} \Het{2}(Y_{\overline{K}},\QQ_\ell)(1).
$$
After passing to a finite extension $L/K$, we may assume that 
${\rm NS}(Y_L)={\rm NS}(Y_{\overline{K}})$.
But then, the $G_L$-representation on $\NS(Y_L)$ is trivial, hence it is also trivial 
on $\Het{2}$, and in particular, unramified.
\qed\medskip

Moreover, the next example shows that also the $G_K$-representation 
on the $\ell$-adic cohomology of the K3 double cover $X$ of an Enriques surface
$Y$ does not detect potential good reduction of $Y$.
This phenomenon is related to flower pot degenerations of Enriques surfaces, 
see \cite[Section 3.3]{Persson} and \cite[Appendix 2]{Persson}.

\begin{Example}
   \label{example: enriques}
  Fix a prime $p\geq5$.
  Consider $\PP^5_{\ZZ_p}$ with coordinates $x_i,y_i$, $i=0,1,2$, and inside it
  the complete intersection of 3 quadrics
  $$
  {\cal X} \,:=\,\left\{\,
  \begin{array}{rrrrrrc}
     & x_1^2 & - x_2^2 & + y_0^2 & & -   y_2^2 & =0\\
     x_0^2 &  & - x_2^2 & & + y_1^2 & -  y_2^2 &=0  \\
     x_0^2 & - e^2 x_1^2 & + x_2^2 & & &- p^2 y_2^2 &=0
   \end{array}
   \right.
  $$
  where $e\in\ZZ_p^\times$ satisfies $e^2\not\equiv0,1,2\mod p$ 
  (for example, we could take $e=2$).
  Then, $\imath:x_i\mapsto x_i, y_i\mapsto -y_i$ defines an involution
  on $\PP^5_{\ZZ_p}$, which induces an involution on ${\cal X}$.
  We denote by $X$ the generic fiber of ${\cal X}$, 
  and by $Y:=X/\imath$ the quotient by the involution.
\end{Example}

\begin{Theorem}
 Let $p\geq5$ and let $X\to Y$ be as in Example \ref{example: enriques}.
 Then, $Y$ is an Enriques surface over $\QQ_p$, such that
 \begin{enumerate}
   \item the K3 double cover $X$ of $Y$ has good reduction, 
   \item the $G_{\QQ_p}$-action on $\Het{2}(X_{\overline{\QQ}_p},\QQ_\ell)$ is 
    unramified for all $\ell\neq p$, 
    \item $Y$ has semi-stable reduction of flower pot type, but
   \item $Y$ does not have potential good reduction.
 \end{enumerate}
\end{Theorem}

\prf
A straightforward computation shows that $X$ is smooth over $\QQ_p$, 
and that $\imath$ acts without fixed points on $X$.
Thus, $X$ is a  K3 surface and $Y$ is an Enriques surface over $\QQ_p$.
The special fiber of ${\cal X}$ is a non-smooth K3 surface with 
4 RDP singularities of type $A_1$ located at $[0:0:0:\pm1:\pm1:1]$.
Then, the blow-up ${{\cal X}}'_1\to{\cal X}$ of the Weil divisor
$\{x_0-ex_1=x_2-py_2=0\}$ defines a simultaneous resolution of the
singularities of ${\cal X}\to\Spec\ZZ_p$, 
and we obtain a smooth model of $X$ over $\ZZ_p$.
In particular, $X$ has good reduction over $\QQ_p$ and 
the $G_{\QQ_p}$-representation on $\Het{2}(X_{\overline{\QQ}_p},\QQ_\ell)$ 
is unramified for all $\ell\neq p$.

Next, let ${{\cal X}}'_2\to{\cal X}$ be the blow-up of the 4 singular points of
the special fiber.
Then, $\imath$ extends to ${{\cal X}}'_2$, and the special fiber
is the union of 4 divisors $E_i$ with the
minimal desingularization $X_p'$ of the special fiber of ${\cal X}$.
The fixed locus of $\imath$ on $X'_p$ is the union of the four
$(-2)$-curves of the resolution.
Moreover, there exist isomorphisms $E_i\cong\PP^1\times\PP^1$
such that $\imath$ acts by interchanging the two factors.
Thus, the quotient ${\cal X}_2'/\imath$ is a model of $Y$
over $\ZZ_p$, whose special fiber is a rational surface
$X_p'/\imath$ (a so-called Coble surface) meeting transversally 
four $\PP^2$'s, that is, a semi-stable degeneration of flower pot
type (see, \cite[Section 3.3]{Persson}).

Seeking a contradiction, we assume that $Y$ has potential good reduction.
Then, there exists a finite extension $L/\QQ_p$ and a smooth model
${\cal Y}\to\Spec\OO_L$ of $Y_L$.
Let ${\cal X}_3\to{\cal Y}\to\Spec\OO_L$ be its 
K3 double cover,
which is a family of smooth K3 surfaces with generic fiber $X_L$, 
whose fixed point free involution specializes to a
fixed point free involution in the special fiber of ${\cal X}_3$.

Now, ${\cal X}_3$ and the base-change of ${\cal X}_1'$ to $\OO_L$
both are smooth models of $X_L$.
The isomorphism of generic fibers extends to a birational
map of special fibers.
The involution on generic fibers extends to rational involutions
of the two special fibers, compatible with the just established
birational map.
Since both special fibers are K3 surfaces, the birational maps
and rational involutions extend to isomorphisms and involutions.
However, in one special fiber the involution acts without fixed points,
whereas it has four fixed curves in the other, a contradiction.
\qed\medskip

\section{Existence and termination of flops}
\label{sec: flops}

Let $X$ be a smooth and proper surface over $K$
with numerically trivial canonical sheaf % $\omega_{X/K}$, 
and assume that we have a smooth model ${\cal X}\to\Spec\OO_K$.
Now, if $\cal L$ is an ample invertible sheaf on $X$, then its specialization 
${\cal L}_0$ to the special fiber may not be ample, and not even be nef.
In this section, we show that there exists a finite sequence of birational
modifications (flops) of ${\cal X}$, such that we eventually 
arrive at a smooth model ${\cal X}^+\to\Spec\OO_K$ of $X$, such that
the restriction of $\cal L$ to the special fiber of ${\cal X}^+$ is big and nef.
We end this section by showing that any two smooth models of $X$
over $\OO_K$ are related by a finite sequence of 
flopping contractions and their inverses.

We start by adjusting \cite[Definition 3.33]{kollar mori} and
\cite[Definition 6.10]{kollar mori} to our situation.

\begin{Definition}
  Let $X$ be a smooth and proper surface over $K$ with numerically 
  trivial $\omega_{X/K}$ that admits a smooth model ${\cal X}\to\Spec\OO_K$.
  Then,
  \begin{enumerate}
    \item a proper and birational morphism $f \colon \cal X \to \cal Y$ over $\OO_K$ 
      is called a {\em flopping contraction}
      if ${\cal Y}$ is normal and if the exceptional locus of $f$ is of codimension at least two.
    \item If $D$ is a Cartier divisor on ${\cal X}$, then a
      birational map ${\cal X}\dashrightarrow{\cal X}^+$ over $\OO_K$
      is called a $D$-{\em flop}
      if it decomposes into a flopping contraction $f \colon \cal X \to \cal Y$ 
      followed by (the inverse of) a flopping contraction $f^+ \colon \cal X^+ \to \cal Y$
      such that $-D$ is $f$-ample and $D^+$ is $f^+$-ample,
      where $D^+$ denotes the strict transform of $D$ on ${\cal X}^+$.
      If ${\cal L}$ is an invertible sheaf on $\cal X$, we
      similarly define an ${\cal L}$-{\em flop}.
   \item A morphism $f^+$ as in (2) is also called a {\em flop} of $f$.
  \end{enumerate} 
\end{Definition}
In general, one also has to assume that $\omega_{\cal X / \OO_K}$
is numerically $f$-trivial in the definition of a flopping contraction.
However, in our situation this is automatic.
Also, a flop of $f$, if exists, does not depend on the choice of $D$ by
\cite[Corollary 6.4]{kollar mori} and \cite[Definition 6.10]{kollar mori}.
This justifies talking about flops without referring to the divisor $D$.

\subsection{Existence of flops}
The following is an adaptation of
Koll\'ar's proof \cite[Proposition 2.2]{Kollar} of the existence of $3$-fold flops
over $\CC$ to our situation, which deals with special 
flops in mixed characteristic.

\begin{Proposition}[Existence of flops]
  \label{prop: existence of flops}
  Let $X$ be a smooth and proper surface over $K$ with numerically 
  trivial $\omega_{X/K}$ that has a smooth model ${\cal X}\to\Spec\OO_K$.
  If ${\cal L}$ is an ample invertible sheaf on $X$ and $C$ is
  an integral (but not necessarily geometrically integral) curve on the special
  fiber ${\cal X}_0$ with ${\cal L}_0\cdot C<0$,
  then there exists a flopping contraction $f \colon {\cal X} \to {\cal X}'$
  and its ${\cal L}$-flop $f^+ \colon {\cal X}^+ \to {\cal X}'$ with the following
  properties
   \begin{enumerate}
      \item $f$ contracts $C$ and no other curves,
      \item ${\cal X}^+\to\Spec\OO_K$ is a smooth model of $X$,
      \item $f$ and $f^+$ induce isomorphisms of generic fibers,
      \item ${\cal L}_0^+\cdot C^+>0$, where ${\cal L}^+$ denotes the extension of ${\cal L}$ on ${\cal X}^+$, and
         where $C^+$ denotes the flopped curve (that is, the exceptional locus of $f^+$).
    \end{enumerate}
\end{Proposition}

\prf
Since ${\cal L}$ is ample, ${\cal L}^{\otimes n}$ is effective for $n\gg0$,
and thus, also its specialization ${\cal L}_0^{\otimes n}$ to the special fiber 
${\cal X}_0$ is effective.
In particular, ${\cal L}_0$ has positive intersection with every ample divisor on ${\cal X}_0$, 
that is, ${\cal L}_0$ is pseudo-effective.
Thus,  there exists a Zariski--Fujita decomposition on $({\cal X}_0)_{\overline{k}}$
$$
  ({\cal L}_0)_{\overline{k}} \,=\, P\,+\,N,
$$
where $P$ is nef, and where $N$ is a sum of effective divisors, 
whose intersection matrix is negative definite, see for example,
\cite[Theorem 14.14]{Badescu}.
Since $\omega_{{\cal X}_0/k}$ is numerically trivial, the adjunction formula shows that
every reduced and irreducible curve in $N$ is a $\PP^1$ with self-intersection $-2$, that is,
a $(-2)$-curve.
Moreover, negative definiteness and the classification of Cartan matrices implies that
$N$ is a disjoint union of ADE curves. 
Next, $k$ is perfect and since the Zariski--Fujita decomposition is unique, it is
stable under $G_k$, and thus, descends to ${\cal X}_0$.

After these preparations, let $C$ be as in the statement, that is, ${\cal L}_0\cdot C<0$.
First, we want to show that there exists a morphism $f:{\cal X}\to{\cal X}'$
of algebraic spaces that contracts $C$.
Being contained in the support of $N$, the base-change
$C_{\overline{k}}\subset({\cal X}_0)_{\overline{k}}$ is a disjoint union of ADE curves. 
Since $C^2<0$, Artin showed that there exists a morphism of projective surfaces over $k$
$$
 f_0\,:\,{\cal X}_0\,\to\,{\cal X}_0'
$$
that contracts $C$ and nothing else (see \cite[Theorem 3.9]{Badescu}, for example).
Since $C_{\overline{k}}$ is a union of ADE-curves, it follows that
$({\cal X}_0')_{\overline{k}}$ has RDP singularities, which are rational and Gorenstein.
Thus, also ${\cal X}_0'$ has rational Gorenstein singularities. 

For all $n\geq0$, we define
$$
  {\cal X}_n \,:=\, {\cal X}\times_{\Spec\OO_K}\Spec( \OO_K/\idealm^{n+1})\,.
$$ 
Since $f_0$ is a contraction with $R^1f_{0 \ast}\OO_{{\cal X}_0}=0$, there exists
a blow-down $f_n:{\cal X}_n\to{\cal X}_n'$ that
extends $f_0$, see \cite[Theorem 3.1]{Cynk van Straten}.
Passing to limits, we obtain a contraction of formal schemes
$$
  \widehat{f}\,:\,\widehat{{\cal X}}\,\to\,\widehat{{\cal X}}'.
$$
By \cite[Theorem 3.1]{artin II}, there exists a contraction of algebraic spaces
$$
 f \,:\, {\cal X}\,\to\,{\cal X}',
$$
whose completion along their special fibers coincides with $\widehat{f}$.
In particular, $f$ is an isomorphism outside $C$ and contracts $C$ to a 
singular point $w\in{\cal X}'$.

Let $\widehat{w}$ be the formal completion of ${\cal X}'$ along $w$, and let
$$
 \widehat{{\cal Z}} \,\to\,\widehat{w}
$$
be the formal fiber over $\widehat{f}$.
Then, $\widehat{w}$ is a formal affine scheme, say $\Spf R$,
and let $k'$ be the residue field, which is a finite extension of $k$.
Let $\OO_{K'}$ be the unramified extension of $\OO_K$ corresponding
to the field extension $k\subseteq k'$.
Since $k\subseteq k'$ is separable, $k'$ arises by adjoining a root $\alpha$
of some monic polynomial $f$ with values $k$.
After lifting $f$ to a polynomial with values in $\OO_K$, 
and using that $R$ is Henselian,
we can lift $\alpha$ to $R$, which shows that $\OO_{K'}$ is contained in $R$.
In particular, we can view $R$ as a local $\OO_{K'}$-algebra without residue
field extension - we denote by $\tilde{R}$ the ring $R$ considered as $\OO_{K'}$-algebra.

Then,
the special fiber of $\Spf\tilde{R}$ is a rational singularity of multiplicity $2$, and thus, by \cite[Lemma 23.4]{Lipman}, 
the completion of the local ring of the special fiber is of the form
\begin{equation}
  \label{eq: RDP equation}
     k'[[x,y,z]]/ \left( h'(x,y,z)\right )\,.
\end{equation}
Using Hensel's lemma, we may assume after a change of coordinates
that the power series $h'(x,y,z)$ is of the form
$z^2-h_1(x,y)z-h_0(x,y)$ for some polynomials $h_0(x,y), h_1(x,y)$.
Using Hensel's lemma again, the completion of  $\tilde{R}$ is of the form
\begin{equation}
 \label{eq: RDP deformation}
    \widehat{\OO}_{K'}[[x,y,z]]/ \left( 
    z^2\,-H_1(x,y)z\,-\,H_0(x,y)
    \right)
\end{equation}
where $H_i(x,y)$ is congruent to $h_i(x,y)$ 
modulo the maximal ideal of  $\widehat{\OO}_{K'}$ for $i=1,2$,
%(here, we have used again $p\neq2$ to obtain (\ref{eq: RDP deformation}) from the
%the versal deformation of (\ref{eq: RDP equation})),
see also \cite[Theorem 4.4]{Kawamata}.
(If $p\neq2$, we may even assume $h_1=0$ and $H_1=0$.)
We denote by $t':\Spf \tilde{R}\to\Spf \tilde{R}$ the involution induced by 
$z\mapsto H_1(x,y)-z$.
It is not difficult to see that $t'$ induces $-{\rm id}$ on local Picard groups, 
see, for example, \cite[Example 2.3]{Kollar}.
Since $R$ is equal to $\tilde{R}$ considered as rings, 
we have established an
involution $t:\widehat{w}\to\widehat{w}$ that induces ${-\rm id}$ on local
Picard groups.
We denote by
$$
  \widehat{{\cal Z}}^+\,\to\,\widehat{w}
$$ 
the composition $t\circ\widehat{f}$.
By \cite[Proposition 2.2]{Kollar}, this gives the desired flop formally.

By \cite[Theorem 3.2]{artin II}, there exists a dilatation $f^+:{\cal X}^+\to{\cal X}'$
of algebraic spaces, such that the formal completion of ${\cal X}^+$ 
along the exceptional locus of $f^+$ is given by the just-constructed
$\widehat{{\cal Z}}^+\to\widehat{{w}}$.
Thus, there exists a birational and rational map
$$
  \varphi\,:\, {\cal X}\,\dashrightarrow\, {\cal X}^+,
$$
which is an isomorphism outside $C$.
From the glueing construction it is clear that
${\cal X}^+$ is a smooth model of $X$ over $\OO_K$.
Finally, from the formal picture above, it is clear that the
restriction of ${\cal L}^+$ to ${\cal X}^+_0$
has positive intersection with the flopped curve $C^+$.
\qed\medskip

\subsection{Termination of flops}
 \label{sec: termination of flops}
Having established the existence of certain flops in mixed characteristic, 
we now show that there is no infinite sequence of them.
To do so, one can adjust the proof of termination of flops from 
\cite[Theorem 6.17 and Corollary 6.19]{kollar mori} over $\CC$ to our situation.
Instead, we give another argument that was kindly  
suggested to us by the referee.

We keep the notations and assumptions of Proposition \ref{prop: existence of flops}.
Then, there are two isomorphisms between the $\ell$-adic cohomology groups
of the special fibers ${\cal X}_0$ and ${\cal X}_0^+$.
\begin{enumerate}
\item The first is by composing the comparison isomorphisms relating the cohomology 
groups of special and generic fibers of $\cal X$ and ${\cal X}^+$
$$
 \alpha\,:\, \Het{2}\left((\cal X_0^+)_{\overline k}, \QQ_\ell\right)(1) \,\cong\, 
   \Het{2}\left(X_{\overline K}, \QQ_\ell\right)(1) 
   \,\cong\, \Het{2}\left((\cal X_0)_{\overline k}, \QQ_\ell\right)(1).
$$
\item Next, the composition $\varphi := (f^+)^{-1} \circ f:{\cal X} \dashrightarrow {\cal X}^+$
induces a birational and rational map of special fibers 
$\varphi_0: {\cal X}_0\dashrightarrow {\cal X}_0^+$,
which extends to an isomorphism, since ${\cal X}_0$ and ${\cal X}_0^+$
are minimal surfaces of Kodaira dimension $\geq0$.
Thus, we obtain a second isomorphism via pullback
$$
 \varphi_0^*\,:\, \Het{2}\left((\cal X_0^+)_{\overline k}, \QQ_\ell\right)(1) 
 \,\cong\,\Het{2}\left((\cal X_0)_{\overline k}, \QQ_\ell\right)(1).
$$
\end{enumerate}
We note that both isomorphisms respect the intersection product 
coming from Poincar\'e duality, that is, they are isometries.
For a $(-2)$-curve $C' \subset (\cal X_0)_{\overline k}$,
we let $[C']$ be the associated cycle class in $\Het{2}((\cal X_0)_{\overline k}, \QQ_\ell)(1)$
and we define the {\em reflection in $C'$} to be the isometry
$$
\begin{array}{ccccc}
r_{C'} &:& \Het{2}\left((\cal X_0)_{\overline k}, \QQ_\ell\right)(1) &\to&  \Het{2}\left((\cal X_0)_{\overline k}, \QQ_\ell\right)(1) \\
&& x &\mapsto& x\,+\, (x\cdot[C'])\,[C']\,.
\end{array}
$$
The following lemma compares the isometries $\alpha$ and $\varphi_0^*$
in terms of reflections in $(-2)$-curves.

\begin{Lemma}
 \label{lem: flop and reflection}
 We keep the notations and assumptions as in Proposition \ref{prop: existence of flops} and 
 denote by $C_1, \ldots, C_m$ the connected components of $C_{\overline k}$.
 Then,
     $$ 
      \alpha \circ (\varphi_0^*)^{-1} \,=\, r_1 \,\cdots\, r_m,
    $$ 
  where
 \begin{enumerate}
   \item either the $C_i$'s are disjoint $(-2)$-curves and $r_i=r_{C_i}$,
   \item or each $C_i$ is the union of two $(-2)$-curves $C_{i,1}$ and $C_{i,2}$ 
    intersecting in one point and
    $r_i=r_{C_{i,1}} r_{C_{i,2}} r_{C_{i,1}} = r_{C_{i,2}} r_{C_{i,1}} r_{C_{i,2}}$.
 \end{enumerate}
\end{Lemma}

\prf
First, we consider the case, where $C$ is geometrically integral.
Let ${\cal Z} \subset {\cal X} \times_{\OO_K} {\cal X}^+$ be the closure 
of the diagonal $\Delta(X)\subset X \times_K X$.
%By \cite[Corollary 2.19]{Saito} (this can be generalized to algebraic space models as in \cite[Section 2.3]{matsumoto2}), 
Then, it is not difficult to see that the isomorphism $\alpha$ is
given by  $x\mapsto \pr_{1,*}([{\cal Z}_0] \cdot \pr_2^*(x))$
(see also Lemma \ref{lem: cycle map} below).
%of the two comparison isomorphisms 
%\[
% \alpha \colon \Het{2}(\cal X_0, \QQ_l)(1) \cong \Het{2}(X_{\overline K}, \QQ_l)(1) \cong \Het{2}(\cal X^+_0, \QQ_l)(1)
%\]
%For the formulation of cycle class map in the context of algebraic spaces see Lemma \ref{lem: cycle map} (**needed??).
We set $\cal U := \cal X \setminus C$ and $\cal U^+ := \cal X^+ \setminus C^+$.
Then, we have a commutative diagram with exact rows
$$
 \xymatrix{
  0 \ar[r] & 
  H^2_{C^+}\left((\cal X_0^+)_{\overline k}, \QQ_\ell\right)(1) \ar[r]\ar[d]^{\alpha'} & 
  \Het{2} \left((\cal X_0^+)_{\overline k}, \QQ_\ell\right)(1) \ar[r]\ar[d]^{\alpha} & 
  \Het{2} \left((\cal U_0^+)_{\overline k}, \QQ_\ell\right)(1) \ar[r]\ar[d]^{\alpha''} & 
  0 \\
  0 \ar[r] & 
  H^2_C \left((\cal X_0)_{\overline k}, \QQ_\ell\right)(1) \ar[r] & 
  \Het{2}\left((\cal X_0)_{\overline k}, \QQ_\ell\right)(1) \ar[r] & 
  \Het{2}\left((\cal U_0)_{\overline k}, \QQ_\ell\right)(1) \ar[r] & 
  0, \\
 }
$$
where $\alpha''$ is also defined by $x\mapsto\pr_{1,*}([\cal Z_0] \cdot \pr_2^*(x))$
and where $\alpha'$ is the map induced by $\alpha$. % and $\alpha''$.
By purity, the left terms are $1$-dimensional and generated by the classes $[C]$ and $[C^+]$, 
respectively.
Moreover, the right terms are canonically isomorphic to the orthogonal complements 
of the left terms.
%Since $\cal Z$ restricted to $(\cal X \setminus C) \times (\cal X^+ \setminus C^+)$ is the diagonal,
Since $\cal Z \rvert_{\cal U \times \cal U^+}$ is the graph of the isomorphism 
$\varphi \rvert_{\cal U} \colon \cal U \to \cal U^+$, it follows that
$\alpha''$ coincides with the pullback by the isomorphism $\varphi_0 \rvert_{\cal U_0} = \varphi \rvert_{\cal U_0} \colon \cal U_0 \to \cal U_0^+$.
Since $\alpha$ is an isometry, so is $\alpha'$, and it maps $[C^+]$ either to $[C]$ or to $-[C]$.
From $\alpha(\cal L^+_0) \cdot \alpha(C^+) = \cal L^+_0 \cdot C^+ > 0$ and 
$\alpha(\cal L^+_0) \cdot C = \cal L_0 \cdot C < 0$, we conclude $\alpha([C^+])=-[C]$.
Putting these observations together, we find $\alpha \circ (\varphi_0^*)^{-1} = r_C$.

Now, we consider the general case.
Since the absolute Galois group $G_k$ acts transitively on the $m$ connected components, 
they are mutually isomorphic.
Since the flops in the disjoint $C_i$ commute, we may assume $m = 1$.
As shown in the proof of Proposition \ref{prop: existence of flops}, 
$C_1 = C_{\overline k}$ is an ADE configuration of $(-2)$-curves.
Since $G_k$ acts on the irreducible components of $C_i$ transitively, it is not difficult to
see from the classification of Dynkin diagrams that only configurations of type
$A_1$ and $A_2$ can occur.
We already treated the $A_1$ case above and thus, we may assume 
an $A_2$-configuration, that is, $C_{\overline k} = C_1 = C_{1,1} \cup C_{1,2}$.
Passing to a finite unramified extension $K'/K$ corresponding to an extension $k'/k$
over which $C$ splits, we consider the following diagram of flops and models 
over $\OO_{K'}$
$$
 \xymatrix{
  \cal X \ar@{-->}[r]^{\varphi_1} \ar@{-->}[rd]_{\varphi_2} & 
  \cal X^{1} \ar@{-->}[r]^{\varphi_2} & 
  \cal X^{12} \ar@{-->}[r]^{\varphi_1} & 
  \cal X^{121} \\
  & 
  \cal X^{2} \ar@{-->}[r]^{\varphi_1} & 
  \cal X^{21} \ar@{-->}[r]^{\varphi_2} & 
  \cal X^{212} , \\
 }
$$
where $\varphi_{j}$ denotes the flop at $C_{1,j}$ or at the corresponding curve on other models
(note that our flops induce isomorphisms between the special fibers).
A straightforward computation shows that $r_{C_{1,1}} r_{C_{1,2}} r_{C_{1,1}}$ 
and $r_{C_{1,2}} r_{C_{1,1}} r_{C_{1,2}}$ both act as $-{\rm id}$ on the one-dimensional subspace 
spanned by $[C_{1,1}] + [C_{1,2}]$ inside $\Het{2}((\cal X_0)_{\overline k}, \QQ_\ell)(1)$ 
and as ${\rm id}$ on its orthogonal complement.
Thus, $\varphi_1 \varphi_2 \varphi_1 \colon \cal X \dashrightarrow \cal X^{121}$ 
and $\varphi_2 \varphi_1 \varphi_2 \colon \cal X \dashrightarrow \cal X^{212}$ both
satisfy the conditions of the flop of the contraction $f$.
Hence, they coincide by the uniqueness of flops, and we set
$\varphi:=\varphi_1\varphi_2\varphi_1=\varphi_2\varphi_1\varphi_2: \cal X \dashrightarrow \cal X^+$.
Clearly, $\varphi$ descends to $\OO_K$ and coincides with the flop in $C$
established in Proposition \ref{prop: existence of flops}.
\qed\medskip

% We can also check that each $\varphi_i$ are as in Proposition \ref{prop: existence of flops},
% i.e. the specialization of $\cal L$ is anti-ample on each contraction.
We define a {\em generalized $(-2)$-curve} on a smooth and proper surface $X$
over a perfect field $k$ to be an integral (but not necessarily geometrically integral)
curve $C\subset X$ such that $C_{\overline{k}}=C_1\cup...\cup C_m$ 
is a disjoint union of ADE curves of type $A_1$ or $A_2$.
We note that $C^2=-2m$, that is, such curves are not necessarily of self-intersection $-2$.
Moreover, we define the reflection $r_C$ in $\NS(X)$ or 
$\Het{2}(X_{\overline{k}},\QQ_\ell)(1)$ to be equal to $r_1\cdots r_m$ as in 
Lemma \ref{lem: flop and reflection}, which is equal to the map
$x\mapsto x + \sum_{i=1}^m (x \cdot [C_i]) [C_i]$.
The following lemma is essentially \cite[Remark 8.2.13]{Huybrechts}.

\begin{Lemma}
 \label{lem: termination of reflections}
Let $X$ be a smooth and projective surface over a perfect field with 
numerically trivial canonical sheaf,
and let $x \in \NS(X)$ be a nonzero effective class with $x^2 \geq 0$.
\begin{enumerate}
\item If $x$ is not nef, then there exists a generalized $(-2)$-curve $C$ with $C \cdot x < 0$
and then, $r_C(x)$ is nonzero and effective.
\item We define a sequence in $\NS(X)$ by setting $x_0:=x$ and if
$x_i$ is not nef, then we choose a generalized $(-2)$-curve $C^i$ with $C^i\cdot x_i<0$ and
set $x_{i+1}:=r_{C^i}(x_i)$.
Then $\{x_i\}$ is a finite sequence of nonzero and effective classes in 
$\NS(X)$ and the last class is nef.
\end{enumerate}
\end{Lemma}

\prf
Using the Zariski--Fujita decomposition of $x$ (see the proof of Proposition \ref{prop: existence of flops})
and Lemma \ref{lem: flop and reflection}, we see that if $x$ is not nef, then 
a generalized $(-2)$-curve $C$ with $C\cdot x<0$ indeed exists.
Since Abelian and bielliptic surfaces do not admit smooth rational curves, 
we may assume that $X$ is a K3 surface or an Enriques surface.
Since $r_C(x)^2 = x^2 \geq 0$, it follows from Riemann--Roch that
either $r_C(x)$ or $-r_C(x)$ is effective.
Let $C_1,...,C_m$ be the connected components of $C_{\overline{k}}$.
Then, we find $x \cdot r_C(x) = x^2 + \sum_i (x \cdot C_i)^2 > 0$,
from which it follows that $x$ and $r_C(x)$ belong to the same 
component of the cone $\{ y \in \NS(X)_{\RR} : y^2 > 0 \}$, and thus, $r_C(x)$ is effective.
This establishes assertion (1).

To show (2), we fix an ample class $H$ of $X$ and let
$C_1,...,C_m$ be the connected components of $C_{\overline{k}}$.
Since $G_k$ acts transitively on these components, we find $x\cdot C_i=\frac{1}{m}x\cdot C$ 
and $H\cdot C_i=\frac{1}{m}H\cdot C$, from which we conclude
$$
r_C(x)\cdot H \,=\, \left(x + \sum_{i=1}^m(x \cdot C_i) C_i\right) \cdot H \,= \,x \cdot H + \frac{1}{m}(x\cdot C)(H\cdot C) \,<\, x \cdot H,
$$
since $x\cdot C<0$ by assumption and $H\cdot C>0$ by ampleness of $H$.
Therefore, if $\{x_i\}\in\NS(X)$ is as in assertion (2), then
$\{x_i \cdot H\}$ is a strictly decreasing sequence of positive integers.
In particular, it must be of finite length, and its last class must be nef.
\qed\medskip

After these preparations, we obtain the following.

\begin{Proposition}[Termination of flops]
  \label{prop: termination of flops}
  Let $(X,{\cal L})$ and ${\cal X}\to\Spec\OO_K$ be as in 
  Proposition \ref{prop: existence of flops}. 
  Then, every sequence of flops as in Proposition \ref{prop: existence of flops}
  is finite.
  In particular, after a finite sequence 
  $$
    ({\cal X},{\cal L}) \,\dashrightarrow\, ({\cal X}^+,{\cal L}^+)
    \,\dashrightarrow\, ({\cal X}^{+2},{\cal L}^{+2})\,\dashrightarrow\,...
    \,\dashrightarrow\, ({\cal X}^{+N},{\cal L}^{+N})
  $$
   of flops we arrive
   at a smooth model $({\cal X}^{+N},{\cal L}^{+N})$ of $X$ over $\OO_K$
   such that the specialization ${\cal L}_0^{+N}$ is big and nef.
\end{Proposition}

\prf
Let $\cdots \dashrightarrow ({\cal X}^{+i},{\cal L}^{+i}) \dashrightarrow \cdots$ 
be a sequence of flops in generalized $(-2)$-curves $C^i\subset ({\cal X}^{+i})_0$ 
as asserted by Proposition \ref{prop: existence of flops} and Lemma \ref{lem: flop and reflection}.
Moreover, we have $[{\cal L}^{+(i+1)}_0] = [r_{C^i}({\cal L}^{+i}_0)]$ by Lemma \ref{lem: flop and reflection}. 
By Lemma \ref{lem: termination of reflections}, this sequence
is finite and ${\cal L}^{+N}_0$ is nef.
\qed\medskip

\subsection{Morphism to a projective scheme}
In the situation of Proposition \ref{prop: termination of flops},
we obtain a birational morphism to a projective scheme as follows.

\begin{Proposition}
  \label{prop: map to a projective scheme}
  Let $X$ be a smooth and proper surface over $K$ with numerically 
  trivial $\omega_{X/K}$ that admits a smooth model ${\cal X}\to\Spec\OO_K$.
  Let ${\cal L}$ be an ample invertible sheaf on $X$ 
  and assume that $\cal L_0$ is big and nef.
  Then, the natural and a priori rational map
  $$
    \pi\,:\, {\cal X}\,\to\,{\cal X}':={\rm Proj}\, \bigoplus_{n\geq0} H^0\left({\cal X},{{\cal L}}^{\otimes n}\right)
  $$
  is a morphism over $\Spec\OO_K$ to a projective scheme.
  More precisely,
  \begin{enumerate}
  \item $\pi$ is a flopping contraction and induces an isomorphism of generic fibers.
  \item The induced morphism on special fibers $\pi_0:{\cal X}_0\to{\cal X}'_0$ 
  is birational and contracts precisely those curves that have zero-intersection with $\cal L_0$.
  In particular, $({\cal X}'_0)_{\overline{k}}$ is a proper surface with at worst
  RDP singularities and $\pi_0$ is the minimal resolution of singularities.
  \end{enumerate}
\end{Proposition}

\prf
Note that also $\omega_{{\cal X}_0/k}$ is numerically trivial.
Since ${\cal L}_0$ is big and nef, we obtain a proper and birational morphism
%(see, for example, \cite[Theorem 14.19]{Badescu})
$$
  \varpi \,:\, {\cal X}_{0} \,\to\, W \,:=\, \Proj \bigoplus_{n\geq0}H^0({\cal X}_{0},{\cal L}_{0}^{\otimes n}).
$$
Base-changing to $({\cal X}_{0})_{\overline{k}}$, the induced morphism $\varpi_{\overline{k}}$
contracts an integral curve $C$ if and only if it has zero-intersection with ${\cal L}_{0}$.
Since the intersection matrix formed by contracted curves is negative definite,
and since an integral curve with negative self-intersection on a surface with numerically
trivial canonical sheaf over an algebraically closed field is a $(-2)$-curve,
%  (that is, $C\iso\PP^1$ and $C^2=-2$), 
it follows from the classification of Cartan matrices,
that $W_{\overline{k}}$ has at worst RDP singularities.

Now, ${\cal L}_{0}^{\otimes n}$ is of degree $0$ on contracted curves for all $n$, and
over $\overline{k}$, these curves are ADE curves.
Thus, we find $R^1\varpi_\ast{\cal L}_0^{\otimes n}=0$ for all $n\geq0$, which implies
$H^1({\cal X}_{0},{\cal L}_0^{\otimes n})=H^1(W,\OO_W(n))$ for all $n\geq0$, and note that
the latter term is zero for $n\gg0$ by Serre vanishing.
Replacing ${\cal L}$ by some sufficiently high tensor power will not change
$\varpi$, and then, we may assume that $H^1({\cal X}_0,{\cal L}_0^{\otimes n})=0$ for all $n\geq1$.
If $f:{\cal X}\to\Spec\OO_K$ denotes the structure morphism, then
semi-continuity and the previous vanishing result imply $R^1f_\ast{\cal L}^{\otimes n}=0$
for all $n\geq1$.
Thus, global sections of ${\cal L}_0^{\otimes n}$ extend to ${\cal L}^{\otimes n}$, and since the former
is globally generated for $n\gg0$, so is the latter.
Thus, we obtain a morphism of algebraic spaces over $\Spec\OO_K$
$$
  \pi\,:\,{\cal X}\,\to\,{\cal X}':={\rm Proj}\, \bigoplus_{n\geq0} H^0\left({\cal X},{{\cal L}}^{\otimes n}\right).
$$
Since ${\cal L}$ is ample on $X$, $\pi$ induces an isomorphism of generic fibers.
Moreover, we can identify the induced map $\pi_0$ on special fibers 
with $\varpi \colon {\cal X}_0\to W$ from above.
\qed\medskip

\subsection{Birational relations among smooth models}
As an application of existence and termination of flops, Koll\'ar \cite[Theorem 4.9]{Kollar}
showed that any two birational complex threefolds with $\QQ$-factorial terminal singularities 
and nef canonical classes are connected by a finite sequence of flops.

We have the following analog in our situation, but since we are dealing with algebraic
spaces rather than projective schemes
(which is analogous to the case of analytic threefolds in loc. cit.), 
the flops as defined above do not suffice.
We only show that two smooth models are connected by flopping contractions and their inverses.

\begin{Proposition}
 \label{prop: flops along models}
 Let $X$ be a smooth and proper surface over $K$ with numerically trivial $\omega_{X/K}$ that has good reduction.
 If ${\cal X}_i\to\Spec\OO_K$ are two smooth models of $X$, then 
 \begin{enumerate}
   \item the special fibers of ${\cal X}_1$ and ${\cal X}_2$ are isomorphic, and
   \item ${\cal X}_1$ and ${\cal X}_2$ are connected by a sequence of birational rational maps 
   that are compositions of flopping contractions and their inverses.
 \end{enumerate}
\end{Proposition}

\prf
The special fibers of ${\cal X}_1$ and ${\cal X}_2$ are birational by the Matsusaka--Mumford theorem \cite[Theorem 2]{MM},
and since they are minimal surfaces of Kodaira dimension $\geq0$, they are isomorphic.
(Note that this statement also follows from the much more detailed analysis below.)

Now, choose an ample invertible sheaf ${\cal L}$ on $X$.
By Proposition \ref{prop: termination of flops}, there exist finite sequences of
flops ${\cal X}_i\dashrightarrow \cdots \dashrightarrow{\cal Y}_i$, $i=1,2$, such that
${\cal L}$ restricts to big and nef invertible sheaves on the special fibers of ${\cal Y}_{i}$.

Applying Proposition \ref{prop: map to a projective scheme}
to our models ${\cal Y}_i$, we obtain flopping contractions
$$
   {\cal Y}_i\,\to\,{\cal Y}_i'\,:=\,\Proj\,\bigoplus_{n\geq0} H^0\left({\cal Y}_i,{\cal L}^{\otimes n}\right)
$$
Now, ${\cal L}$ is ample on ${\cal Y}_i'$. 
Moreover, the ${\cal Y}_i'$ are normal projective schemes and birational outside a finite 
number of curves in their special fibers.
In fact, there exists a birational and rational map between them that is compatible with ${\cal L}$.
Thus, by \cite[Theorem 5.14]{Kovacs}, this birational map extends to an isomorphism,
and then, we obtain a birational map ${\cal Y}_1\dashrightarrow {\cal Y}_2$
with decomposition ${\cal Y}_1 \to {\cal Y}_1' \iso {\cal Y}_2' \leftarrow {\cal Y}_2$ 
of the required form.

Putting all these birational modifications together, we have connected 
${\cal X}_1$ and ${\cal X}_2$ by a sequence of birational maps of the required form.
\qed\medskip

\section{Group actions on models}
\label{sec: group actions on models}

In this section, we study group actions on models.
More precisely, we are given a smooth and proper
surface $X$ over $K$ with numerically trivial $\omega_{X/K}$,
a finite field extension $L/K$,
which is Galois with group $G$, 
and a smooth proper model ${\cal X}\to\Spec\OO_L$ of $X_L$.
Then we study the following questions:
\begin{enumerate}
 \item Does the $G$-action on $X_L$ extend to ${\cal X}$? 
 \item If so, is the special fiber $({\cal X}/G)_0$ of the quotient equal to 
    the quotient ${\cal X}_0/G$ of the special fiber? 
\end{enumerate}
It turns out, that the answer to question (1) is ``yes'',
when allowing certain birational modifications of the model,
and in question (2), it turns out that the case where $p\neq0$ and $p$ 
divides the order of $G$ (wild group actions) is subtle.

\subsection{Extending group actions to a possibly singular model}
Given a smooth and proper surface $X$ over $K$ with numerically trivial $\omega_{X/K}$
that admits a model ${\cal X}$ with good reduction
after a finite Galois extension $L/K$ with group $G$, we
first show that the $G$-action extends to a (mild) birational
modification of ${\cal X}$.

\begin{Proposition}
 \label{prop: rational group action}
  Let $X$ be a smooth and proper surface over $K$
  with numerically trivial $\omega_{X/K}$.
  Assume that there exist 
  \begin{enumerate}
   \item a finite Galois extension $L/K$ with Galois group $G$, as well as 
   \item a smooth model ${\cal X}\to\Spec\OO_L$ of $X_L$, and
   \item an ample invertible sheaf ${\cal L}$ on $X$,
       whose pull-back to $X_L$ restricts to an invertible sheaf on the special
       fiber ${\cal X}_0$ that is big and nef.
  \end{enumerate}
  Then, there exists a proper birational morphism $\pi$
  $$
     \xymatrixcolsep{3pc}
     \xymatrix{
       {\cal X} \ar[d]\ar[r]^{\pi} & \ar[dl]{\cal X}' \\
       \Spec\OO_L }
   $$
  of algebraic spaces over $\OO_L$, such that
  \begin{enumerate}
    \item The natural $G$-action on $X_L$ 
    extends to ${\cal X}'$ and is compatible
    with the $G$-action on $\OO_L$.
    \item ${\cal X}'$ is a projective scheme over $\Spec \OO_L$.
    \item The generic fibers of ${\cal X}$ and ${\cal X}'$ are isomorphic via $\pi$, whereas the 
       induced morphism on special fibers
       $\pi_0:{\cal X}_0\to{\cal X}_0'$
       is birational and projective, such that the geometric
       special fiber $({\cal X}_0')_{\overline{k}}$ has at worst RDP singularities.
%    \item If $\pi$ is not an isomorphism, then ${\cal X}'$ is not regular.
  \end{enumerate}
Moreover, if $\pi$ is not an isomorphism, then ${\cal X}'$ is not regular.
\end{Proposition}

\proof
Since ${\cal X}$ is regular, the pull-back of ${\cal L}$ to $X_L$ extends
to an invertible sheaf on ${\cal X}$.
By abuse of notation, we shall denote the pull-back to $X_L$ and its extension
to ${\cal X}$ again by ${\cal L}$.
By assumption, the restriction ${\cal L}_0$ of ${\cal L}$ to the special fiber
${\cal X}_0$ is big and nef.
Note that $\omega_{{\cal X}_0/k}$ is numerically trivial.
Let 
$$
  \pi\,:\,{\cal X}\,\to\,{\cal X}':={\rm Proj}\, \bigoplus_{n\geq0} H^0({\cal X},{{\cal L}}^{\otimes n})
$$
be the morphism of algebraic spaces over $\OO_L$ given by Proposition \ref{prop: map to a projective scheme}.
Then, $\pi$ has all the properties asserted in claim (3) of the proposition.
Clearly, ${\cal X}'$ is a projective scheme over $\Spec\OO_L$, and if
$\pi$ is not an isomorphism, then the exceptional locus is non-empty
and of codimension $2$, which implies that ${\cal X}'$ cannot be
regular by van~der~Waerden purity, see, for example,
\cite[Theorem 7.2.22]{Liu}.

It remains to establish the $G$-action on ${\cal X}'$.
Since ${\cal L}$ is a $G$-invariant invertible sheaf on $X_L$,
we have an induced $G$-action on $H^0(X_L,{\cal L}^{\otimes n})$ for
all $n\geq0$.
We will show that this extends to an action on $H^0(\cal X,{\cal L}^{\otimes n})$.
%%Let $U\subset{\cal X}$ be an open and dense subspace such that the $G$-action 
%%extends to $U$.
%Since birational rational maps between normal spaces are isomorphisms outside subspaces of codimension $\geq 2$
%(use, for example, \cite[Corollary 4.4.3]{Liu})
%and since $G$ is finite, 
%%Since ${\cal X}$ is normal, we may assume that the complement 
%we can take an open and dense subspace $U\subset{\cal X}$ such that the $G$-action 
%extends to $U$, 
%the complement ${\cal X}\backslash U$ is of codimension $\geq2$, and ${\cal X}\backslash U$ is contained in ${\cal X}_0$.
First, we show that there exists a closed subspace ${\cal Z}\subset\cal X$ of codimension $\geq 2$ that is
contained in $\cal X_0$, such that the $G$-action on $X_L$ 
extends to an action on ${\cal U} := {\cal X} \setminus {\cal Z}$.
Since every birational rational map between two normal algebraic spaces 
is an isomorphism outside a closed subspace of codimension $\geq 2$, % (use, for example, \cite[Corollary 4.4.3]{Liu}),
there exists for every $g\in G$ a closed subspace ${\cal Z}_g\subset{\cal X}$ of codimension $\geq 2$ that is 
contained in $\cal X_0$ and such that $g:X_L \to X_L$ extends to $g \colon {\cal X} \setminus{\cal Z}_g \to \cal X$.
Since $\cal X_0$ is a minimal surface, the restriction 
$g \rvert_{{\cal X}_0 \setminus {\cal Z}_g}:{\cal X}_0 \setminus {\cal Z}_g \to {\cal X}_0$ to the special fiber
extends to an automorphism $g: {\cal X}_0 \to{\cal X}_0$.
This defines a $G$-action on ${\cal X}_0$.
Let ${\cal Z}' := \bigcup_{g \in G} {\cal Z}_g$ and ${\cal Z} := \bigcup_{g \in G} g({\cal Z}')$, 
where $g({\cal Z}')$ is the image by the action just defined.
Since $G$ is a finite group, ${\cal Z}$ is closed.
This ${\cal Z}$ satisfies the above condition, for otherwise there exists a $g \in G$,
such that the image of $g:{\cal X}\setminus {\cal Z} \to \cal X$ is not contained in 
${\cal X} \setminus {\cal Z}$, or, equivalently $g^{-1}({\cal Z}) \not\subset {\cal Z}$.
%a locally-closed subspace $W \subset Z$ such that $g^{-1}
%a (possibly non-proper) curve $C \subset U$ such that $g(C) \not\subset U$,
%or equivalently $g(C) \subset Z$,
However, since ${\cal Z}$ is $G$-stable, this cannot happen.
Thus we obtain a $G$-action on ${\cal U} = {\cal X} \setminus {\cal Z}$.

If $s$ is a global section of ${\cal L}^{\otimes n}$ over ${\cal X}$ and $\sigma\in G$,
then $\sigma(s|_{\cal U})$ is a well-defined global section of ${\cal L}^{\otimes n}$ over 
${\cal U} = \cal X \setminus \cal Z$.
Since ${\cal L}^{\otimes n}$ is a reflexive sheaf on a regular algebraic space, 
$\sigma(s|_{\cal U})$ extends uniquely to a global section of ${\cal L}^{\otimes n}$ over $\cal X$.
Thus, we obtain a $G$-action on $H^0({\cal X},{\cal L}^{\otimes n})$, which gives rise to 
a $G$-action on ${\cal X}'$ that is compatible with the $G$-action on $\OO_L$, as well as with
the  natural $G$-action on $X_L$.
\qed\medskip

\begin{Remark}
  If all assumptions of Proposition \ref{prop: rational group action}
  except assumption (3) are satisfied, then Proposition \ref{prop: termination of flops}
  shows that there exists another smooth model of $X$ over $\Spec \OO_L$
  for which all assumptions including (3) hold, 
  and to which we can apply Proposition \ref{prop: rational group action}.
\end{Remark}

\subsection{Examples where the action does not extend}
In general, it is too much to ask for an extension of the $G$-action 
from $X_L$ to ${\cal X}$ (notation as in Proposition \ref{prop: rational group action}).
The following example is typical.

\begin{Example}[Arithmetic $3$-fold flop]
 Consider $\QQ_p$ with $p\neq2$ and set $L:=\QQ_p(\varpi)$, where
 $\varpi^2=p$.
 Then, $L/\QQ_p$ is Galois with group $G=\ZZ/2\ZZ$
 and the non-trivial element of $G$ acts as $\varpi\mapsto-\varpi$.
 We equip
 $$
   {\cal X}'\,:=\,\Spec \OO_L[[x,y,z]]/(xy+z^2-\varpi^2)\,\to\,\Spec\OO_L
 $$
 with the $G$-action that is the Galois action on $\OO_L$, and that is
 trivial on $x,y,z$.
 It is easy to see that the induced $G$-action on the special fiber ${\cal X}'_0$
 is trivial.
 Next, we consider the two ideal sheaves ${\cal I}_{\pm}:=(x,y,z\pm\varpi)$ 
 of $\OO_{{\cal X}'}$ and their blow-ups
 $$
   \pi_{\pm}\,:\,{\cal X}_\pm\,\longrightarrow\,{\cal X}'
 $$
 Then, ${\cal X}_\pm$ are regular schemes, ${\cal X}'$
 is singular at the closed point $(x,y,z,\varpi)$, $\pi_{\pm}$
 are both resolutions of singularities, and the exceptional locus 
 is a $\PP^1$ in both cases.
 The ideals ${\cal I}_\pm$ are not $G$-invariant and the
 $G$-action on ${\cal X}'$ does not extend to that on ${\cal X}_+$ nor on ${\cal X}_-$.
 (Instead, the non-trivial element of $G$ induces an isomorphism 
 ${\cal X}_+ \to {\cal X}_-$.)
 In fact, ${\cal X}'$ is an arithmetic version of a $3$-fold ordinary 
 double point, and the rational map ${\cal X}_+\dashrightarrow{\cal X}_-$
 is an arithmetic version of the classical Atiyah flop.
\end{Example}

Even worse, the following example (which is a modification of 
Example \ref{example: no good reduction} below, 
and rests on examples from \cite[Section 5.3]{matsumoto2} and \cite[Section 3]{van luijk})
shows that if we have a $G$-action
on a singular model ${\cal X}'$ as in Proposition \ref{prop: rational group action},
then there may exist resolutions of singularities to which the $G$-action extends,
as well as resolutions to which the $G$-action does not extend.
Moreover, our examples are models of K3 surfaces, that is, such phenomena
are highly relevant for our discussion.

\begin{Example}
\label{example: resolutions} 
    Before giving explicit examples, let us explain the strategy,
    again, for $K=\QQ_p$ and $k=\FF_p$ in this example:
    
    Assume $p\neq2$, let $k'/k$ be the unique extension of degree $2$,
    let $K'/K$ be the corresponding unramified extension, and let
    $G=\{1,\sigma\}$ be the Galois group of both extensions.
    Next, let ${\cal X}' \to \Spec \OO_K$ be a proper scheme such that 
    the geometric special fiber $({\cal X}'_0)_{\overline{k}}$ has only RDP singularities 
    and at least two of them. 
    Set $S:=\Sing {\cal X}'_0$, and assume that all points of $S$ are $k'$-rational but not all 
    $k$-rational.
    Then, $G$ acts non-trivially on $S(k')= S(\overline{k})$.
     Let us finally assume that there exist two different resolutions of singularities
     $\psi_{\pm} \colon {\cal X}_{\pm} \to {\cal X}'$, both of which are isomorphisms 
     outside $S$, and both of which are
     obtained by blowing up ideal sheaves ${\cal I}_{\pm}$ defined over $\OO_K$.
     From this setup, we can produce the announced counterexamples:
     \begin{enumerate}
       \item The  Galois action on $({\cal X}')_{\OO_{K'}}$ extends to $({\cal X}_+)_{\OO_{K'}}$, as well as to
          $({\cal X}_-)_{\OO_{K'}}$. 
          Thus, there do exist resolutions of singularities to which the $G$-action extends.
       \item On the other hand, for each decomposition $S(k') = S_1 \sqcup S_2$,
         we define
         $\psi_{S_1,S_2}:{\cal X}_{S_1,S_2} \to ({\cal X}')_{\OO_{K'}}$ to be the morphism
         that is equal to $\psi_+$ (resp.\ $\psi_-$) on $({\cal X}')_{\OO_{K'}} \setminus S_2$ (resp.\ on $({\cal X}')_{\OO_{K'}} \setminus S_1$).
         We note that $\psi_{S_1,S_2}$ is also a resolution of singularities.
         But now, if $S_1$ and $S_2$ are not $G$-stable, 
         then the $G$-action on $({\cal X}')_{\OO_{K'}}$ does not extend to ${\cal X}_{S_1,S_2}$,
          but induces an isomorphism from ${\cal X}_{S_1,S_2}$ to ${\cal X}_{\sigma(S_1),\sigma(S_2)}$,
          where $\sigma\in G$ is the non-trivial element.
   \end{enumerate}
   
   We now give explicit examples for $p\geq5$ and $K=\QQ_p$.
   Fix a prime $p\geq5$ and choose an integer $d$ such that $d$ is not a quadratic residue modulo $p$,
   and such that $d^6 \not\equiv -2^{-4}\cdot 3^{-3} \mod p$ (one easily checks that such $d$ exists). 
   We define the polynomial
   $$
     \begin{array}{ccl}
       \phi &:=& x^3 -x^2y -x^2z+x^2w-xy^2-xyz+2xyw+xz^2+2xzw \\ 
         && {}+ y^3+y^2z-y^2w+yz^2+yzw-yw^2+z^2w+zw^2+2w^3.
     \end{array}
    $$
    Then, we choose a homogeneous polynomial $f\in\ZZ[x,y,z,w]$
    of degree $3$, such that the following congruences hold
    $$
      \begin{array}{ccll}
         f &\equiv& \phi & \mod 2\\
         f &\equiv& z(x^2-z^2)+w^3 &\mod p.
      \end{array}
    $$
    Next, we choose homogeneous quadratic polynomials
    $2g, 2h \in \ZZ[x,y,z,w]$, such that the following congruences hold
    $$
       \begin{array}{ccll}
          g^2 - p^2 h^2 &\equiv& (z^2 + xy + yz) (z^2 + xy) & \mod 2 \\
          g^2 - p^2 h^2 &\equiv& (y^2 - d x^2)^2  & \mod p.
       \end{array}
    $$
     Finally, we define the quartic hypersurface
     $$
       {\cal X}'\,:=\,{\cal X}'(p) \,:=\, \left\{ wf + g^2 - p^2 h^2 = 0
      \right\} \,\subset\, \PP^3_{\ZZ_p},
    $$
    and denote by $X=X(p)$ its generic fiber.
    Then, $X$ is a smooth K3 surface over $\QQ_p$.
% (this follows from the congruences modulo $2$, as in \cite[Section 5.3]{matsumoto2}).
    The subscheme $S = \Sing {\cal X}'_0$ is given by 
   $$
      S \,=\, \left\{ w = y^2 - dx^2 = z(x^2-z^2) = 0 \right\} \,\subset\,\PP^3_{\FF_p}.
   $$ 
   Thus, we find 6 RDP singularities on $({\cal X}_0')_{\overline{\FF}_p}$,
   all of which are defined over $\FF_{p^2}=\FF_p[\sqrt{d}]$, and 
   $G=\Gal(\FF_{p^2}/\FF_p)$ acts non-trivially on $S(k')$, since $\sqrt{d} \not\in\FF_p$.
   Finally, the blow-ups $\psi_{\pm} \colon {\cal X}_{\pm} \to {\cal X}'$
   of the ideals ${\cal I}_{\pm} := (w = g \pm p h = 0)$ 
   are both resolutions of singularities.
   As explained in the strategy above, this setup yields the desired examples.
   (We refer to Remark \ref{rem:phi} for the reason why we use this $\phi$.)
\end{Example}

\subsection{Extending the inertia action to the smooth model}
Despite all these discouraging examples, there are situations, in which
the $G$-action on $X_L$ {\em does} extend to ${\cal X}$,
and not merely to a singular model ${\cal X}'$
(notation as in Proposition \ref{prop: rational group action}).
More precisely, we have the following result.

\begin{Proposition}
  \label{prop: extending group action to smooth model}
  We keep the notations and assumptions of Proposition \ref{prop: rational group action}.
%  We denote by $I_G$ the inertia subgroup of $G$.
  \begin{enumerate}
    \item If $X$ is an Abelian surface or a hyperelliptic surface, then the $G$-action
       on $X_L$ extends to ${\cal X}$.
    \item Let $H \subseteq G$ be a subgroup, whose action on $\Het{2}(X_{\overline L},\QQ_\ell)$ is trivial
      (for example, this is the case if $\Het{2}(X_{\overline L},\QQ_\ell)$ is unramified and $H \subseteq I_G$).
      Then, the $H$-action on $X_L$ extends to ${\cal X}$.
%    If the $I_G$-action on $\Het{2}(X_{\overline L},\QQ_\ell)$ is trivial, that is, the $G$-action is unramified,
%       then the $I_G$-action on $X_L$ extends to ${\cal X}$.
  \end{enumerate}
\end{Proposition}

%We first prove the following lemma.
%We first introduce the following generalization of cycle class maps to the context of algebraic spaces.
We first introduce cycle class maps in the context of algebraic spaces.

\begin{Lemma} 
 \label{lem: cycle map}
 Let ${\cal X}$ be a proper and smooth algebraic space over the spectrum
 $S = \Spec \OO_K$ of a strictly Henselian DVR $\OO_K$, 
 and let ${\cal Z} \subset {\cal X}$ be a closed subspace of codimension $c$ 
  that is flat over $S$.
  Then, the natural isomorphism 
  $$
    \Het{2c}\left(({\cal X}_0)_{\overline k},\, \ZZ/n\ZZ(c)\right) \,\to\, 
    \Het{2c}\left({\cal X}_{\overline K},\, \ZZ/n\ZZ(c)\right)
 $$
 maps $[{\cal Z}_0]$ to $[{\cal Z}_K]$.
\end{Lemma}

\prf
This is immediate if we define a cycle class 
$[{\cal Z}]$ in $\Het{2c}({\cal X}, \ZZ/n\ZZ(c))$ in such a way, that it is
compatible with base change $S' \to S$ of base schemes 
(and thus, in particular, compatible with restrictions to generic and special fibers).
If ${\cal X}$ is a scheme, this is defined and shown in \cite[Cycle, Num\'ero 2.3]{SGA4 1/2},
and thus, it remains to treat the case, where ${\cal X}$ is an algebraic space.

First, let us recall cohomological descent.
Let $V \to {\cal Y}$ be an \'etale covering of schemes.
Such a morphism is of cohomological descent by \cite[Proposition V$^\mathrm{bis}$.4.3.3]{SGA4},
and we have a spectral sequence \cite[Proposition V$^\mathrm{bis}$.2.5.5]{SGA4}
$$
 E_1^{p,q}  \,:=\, \Het{q}(V_p, a_p^* F) \,\Rightarrow\, \Het{p+q}({\cal Y},F),
$$
where, for each $p \geq 0$, $V_p$ is the $(p+1)$-fold fibered product of $V$ over ${\cal Y}$, 
and where $a_p:V_p \to {\cal Y}$ is the structure map.
Next, we consider the more general case, where $V \to {\cal Y}$ is an \'etale covering of an algebraic space by a scheme.
%Now assume $Y$ is an algebraic space but not necessarily a scheme 
%(we assume $V$, and hence $V_p$, is a scheme).
%$V \to Y$ is a \'etale covering of an algebraic space by a scheme
(Note that then, the $V_p$ are schemes.)
We observe that $V \to {\cal Y}$ is still of cohomological descent (proved in the same way) and we obtain the same spectral sequence.

After these preliminary remarks, let $U\to {\cal X}$ be an \'etale covering  by a scheme $U$.
Applying the previous paragraph to the covering $U \times_{\cal X} {\cal Z} \to {\cal Z}$ and the sheaf $R^{2c} i^! \ZZ/n\ZZ(c)$ on ${\cal Z}$ 
(where $i \colon {\cal Z} \to {\cal X}$ is the closed immersion),
and using the isomorphism $H^{2c}_{\cal Z}({\cal X}, \ZZ/n\ZZ(c)) \cong H^0({\cal Z}, R^{2c} i^! \ZZ/n\ZZ(c))$ 
(this isomorphism is also proved by reducing to the scheme case),
we obtain an isomorphism
$$
 \Ker \left( H^{2c}_{Z_0}(U_0, \ZZ/n\ZZ(c)) \,\to\, H^{2c}_{Z_1}(U_1, \ZZ/n\ZZ(c)) \right) 
 \,\cong\, H^{2c}_{\cal Z}({\cal X}, \ZZ/n\ZZ(c)),
$$
where $Z_p = U_p \times_{\cal X} {\cal Z}$ (note that this is a scheme).
We define $[{\cal Z}] \in H^{2c}_{\cal Z}({\cal X}, \ZZ/n\ZZ(c))$ to be $[Z_0] \in H^{2c}_{Z_0}(U_0, \ZZ/n\ZZ(c))$ (this class lies indeed in the kernel).
Since the cycle map (for schemes) is \'etale local, this construction does not depend on the choice of the \'etale covering $U \to {\cal X}$.

Compatibility with change of base schemes reduces to the scheme case.
\qed\medskip

\prf[Proof of Proposition \ref{prop: extending group action to smooth model}]
(1)
It follows from the assumptions that ${\cal X}_0$ is a smooth and proper surface with numerically trivial
$\omega_{{\cal X}_0/k}$, and that it has the same $\ell$-adic Betti numbers as $X$.
Thus, by the classification of surfaces (see, for example, \cite[Theorem 6 and the following Proposition]{bm2}), 
also ${\cal X}_0$ is Abelian and (quasi-)hyperelliptic, respectively.
As seen in the proof of Proposition \ref{prop: rational group action}, the
(geometric) exceptional locus of $\pi$ is a union of $\PP^1$'s with self-intersection number $(-2)$.
Now, there are no rational curves on Abelian varieties.
Also, it follows from the explicit classification and description 
of (quasi-)hyperelliptic surfaces in \cite[Proposition 5]{bm2} that 
they do not contain any smooth rational curves. 
In particular, $\pi$ must be an isomorphism, which implies that 
the $G$-action extends to ${\cal X}$.

(2)
%After replacing $G$ with $I_G$, we may assume that $L/K$ is totally ramified.
After replacing $K$ by an intermediate extension, we may assume $H = G$.
To show that the $G$-action extends, it suffices to show that 
the $\sigma$-action extends for every $\sigma\in G$.
Thus, let $\sigma\in G$, and after replacing $G$ by the cyclic subgroup
generated by $\sigma$, we may assume that $G$ is cyclic, 
say $G=\Gal(L/K)\iso\ZZ/n\ZZ$, and generated by $\sigma$.

Let $U\subset{\cal X}$ be the maximal open subspace to which the $G$-action on $X_L$ extends.
Then, as in the proof of Proposition \ref{prop: rational group action}, $U$ contains the generic fiber $X_L$, as well as an open dense subscheme 
of the special fiber ${\cal X}_0$.
Let $\Gamma \subset \X^n$ be the closure of the set 
$\{ (x, \sigma(x), \ldots, \sigma^{n-1}(x)) \mid x \in U \} $
in $\X^n$.
The group $G = \ZZ/n\ZZ$ acts on $\X^n$ by permutation of the factors,
and this action restricted to $\Gamma_L$ coincides with the natural $G$-action on $X_L$
via $\pr_1 \colon \Gamma_L \stackrel{\sim}{\to} X_L$.

Consider the diagram of $G$-representations
\begin{equation}
 \label{eq: cycle diagram}
\xymatrix{
    H^2({\cal X}_{\overline 0}^n) \ar[dd]^{\cdot [\Gamma_{\overline 0}]} & H^2(X_{\overline L}^n) \ar[l]_\iso \ar[rd] \ar[dd]^{\cdot [\Gamma_{\overline L}]} & \\
& & H^2(\Gamma_{\overline L}) \iso H^2(X_{\overline L}) \ar[ld] \\
H^{4n-2}({\cal X}_{\overline 0}^n)  & H^{4n-2}(X_{\overline L}^n), \ar[l]_\iso & \\
}
\end{equation}
where we omit the coefficients (Tate twists of $\QQ_\ell$) of the $\ell$-adic cohomology groups 
from the notation, 
and we write $\Gamma_{\overline 0} := (\Gamma_0)_{\overline k}$ and so on.
The triangle on the right is clearly commutative.
The commutativity of the square on the left 
follows from the fact that the classes $[\Gamma_{\overline 0}]$ and $[\Gamma_{\overline L}]$ correspond 
via the isomorphism $H^{4n-4}({\cal X}_{\overline 0}^n) \cong H^{4n-4}(X_{\overline L}^n)$, 
as proved in Lemma \ref{lem: cycle map}.
%which will be proved in Lemma \ref{lem: cycle map} below.
Note that all irreducible components of $\Gamma_{\overline 0}$ are of dimension $2$
(by, for example, \cite[Proposition 4.4.16]{Liu}).

Let $\pi:{\cal X}\to{\cal X}'$ be as in Proposition \ref{prop: rational group action}.
Let $E \subset {\cal X}_0$ be the exceptional locus of 
$\pi_0:{\cal X}_0 \to {\cal X}'_0$
and $E^\alpha$ be the irreducible components of $E_{\overline k}$.
Each irreducible component $E^\alpha$ is isomorphic to $\PP^1$.
Since $\pi_0$ is a resolution of singularities, the intersection
matrix $(E^\alpha \cdot E^\beta)_{\alpha,\beta}$ is 
negative definite (by Hodge index theorem), hence invertible.
% it is a block matrix, all of whose blocks are negative definite,
In particular, if we are given $c_\alpha \in \QQ_\ell$ for $\alpha=1,...,m$, such that 
$\sum_{\alpha=1}^m c_\alpha E^\alpha \cdot E^\beta = 0$ for all $\beta$, then 
$c_\alpha = 0$ for all $\alpha=1,...,m$.

Consider the irreducible components of $\Gamma_{\overline 0}$.
First, there is the ``diagonal'' component, that is, 
the closure of the set $\{ (x, \sigma(x), \ldots, \sigma^{n-1}(x)) \mid x \in U_0 \}$.
If $Z$ is a non-diagonal component (assuming there is one), then
$Z$ is contained in $E^{\alpha_1} \times \cdots \times E^{\alpha_n}$
for some $\alpha_1, \ldots, \alpha_n$.
From the K\"unneth formula and the fact that $H^*(\PP^1, \QQ_\ell) \iso \QQ_\ell [\PP^1] \oplus \QQ_\ell [\pt]$, 
it follows that the cycle class $[Z] \in H^{4n-4}({\cal X}_{\overline 0}^n)$ is a non-zero 
$\ZZ_{\geq 0}$-combination of $[E^\gamma_i \times E^\delta_j \times \pt^{n-2}]$ with $i \neq j$,
where 
$$
E^\gamma_i \times E^\delta_j \times \pt^{n-2} \,:=\,
\cdots \times E^\gamma \times \cdots \times E^\delta \times \cdots 
$$
(the $i$.th component is equal to $E^\gamma$, the $j$.th component is equal to $E^\delta$, and the remaining components are 
equal to a point).
Hence, if we set $[\Gamma_{\overline 0}]^{\nondiag}:=[\Gamma_{\overline 0}]-[\diag]$, then there exist
$c_{i,j,\gamma,\delta}\in\ZZ_{\geq0}$ such that
\begin{equation}
 \label{eq: nondiag cycle}
  [\Gamma_{\overline 0}]^{\nondiag} \,=\, 
  % [\Gamma_{\overline 0}] - [\diag] \,=\, 
  \sum_{i,j,\gamma,\delta} c_{i,j,\gamma,\delta} \, [E^\gamma_i \times E^\delta_j \times \pt^{n-2}]
  \,\in\,H^{4n-4}({\cal X}_{\overline 0}^n).
\end{equation}
We have $c_{i,i,\gamma,\delta} = 0$ for all $i,\gamma,\delta$.

We want to show $[\Gamma_{\overline 0}]^{\nondiag}=0$.
For this, we will use the assumption that the $G$-action on $H^2(X_{\overline L})$ is trivial.
Using the commutative diagram (\ref{eq: cycle diagram}),  we see that the map
$\cdot [\Gamma_{\overline 0}] : H^{2}({\cal X}_{\overline 0}^n) \to H^{4n-2}({\cal X}_{\overline 0}^n)$
factors through $H^2(X_{\overline L})$, and thus,
every element in its image is $G$-invariant.
In particular, for all $\alpha$ and $i$, the cycle
$[E^\alpha_i \times {\cal X}_{\overline 0}^{n-1}] \cdot [\Gamma_{\overline 0}] \in H^{4n-2}({\cal X}_{\overline 0}^n)$ is $G$-invariant,
where 
$$
E^\alpha_i \times {\cal X}_{\overline 0}^{n-1} \,:=\,
\cdots \times E^\alpha \times \cdots 
$$
(the $i$.th component is equal to $E^\alpha$ and the remaining components are equal to ${\cal X}_{\overline 0}$).
Now, $G$ acts by $\sigma \colon [E^\beta_j \times {\cal X}_{\overline 0}^{n-1}] \mapsto [E^\beta_{j+1} \times {\cal X}_{\overline 0}^{n-1}]$, which implies
that for all $\beta$, the cycle 
$[E^\alpha_i \times {\cal X}_{\overline 0}^{n-1}] \cdot [\Gamma_{\overline 0}] \cdot [E^\beta_j \times {\cal X}_{\overline 0}^{n-1}] \in H^{4n}({\cal X}_{\overline 0}^n) \cong \QQ_\ell$
is independent of $j$.
Since $[E^\alpha_i \times {\cal X}_{\overline 0}^{n-1}] \cdot [\diag] \cdot [E^\beta_j \times {\cal X}_{\overline 0}^{n-1}]$
is equal to $E^\alpha\cdot E^\beta$, it is also independent of $j$,
and thus, $[E^\alpha_i \times {\cal X}_{\overline 0}^{n-1}] \cdot [\Gamma_{\overline 0}]^{\nondiag} \cdot [E^\beta_j \times {\cal X}_{\overline 0}^{n-1}]$
is independent of $j$ for all $\beta$.
In order to compute its value, % (via the isomorphism $H^{4n}({\cal X}_0^n) \cong \QQ_\ell$)
we use equation (\ref{eq: nondiag cycle}) and find 
$$
  [E^\alpha_i \times {\cal X}_{\overline 0}^{n-1}] \cdot [\Gamma_{\overline 0}]^\nondiag \cdot [E^\beta_j \times {\cal X}_{\overline 0}^{n-1}] 
= \sum_{\gamma,\delta} (c_{i,j,\gamma,\delta} + c_{j,i,\delta,\gamma}) 
(E^\gamma \cdot E^\alpha) (E^\delta \cdot E^\beta) .
$$
Since $c_{i,i,\gamma,\delta}=0$ for all $i$, this sum is zero for $i = j$. 
Since it is independent of $j$, this sum is zero for all $i,j$.
Using invertibility of the matrix $(E^\alpha \cdot E^\beta)$ twice, 
we obtain $c_{i,j,\gamma,\delta} + c_{j,i,\delta,\gamma} = 0$ for all $i,j,\gamma,\delta$.
Thus, $[\Gamma_{\overline 0}]^{\nondiag} = 0$.

Now, $\pr_i:\Gamma\to{\cal X}$ is a proper birational morphism for all $i$, where $\cal X$ is regular,
and $\Gamma$ is integral.
Thus, by van der Waerden purity (see, \cite[Theorem 7.2.22]{Liu}, for example, and note that this result
can easily be extended to algebraic spaces),
the exceptional locus of $\pr_i$ is either empty or a divisor.
If it was a divisor, it would give rise to a non-diagonal component of $\Gamma_{\overline{0}}$, 
which does not exist by the previous computations.
Thus, $\pr_i$ is an isomorphism for all $i$, and since the $\Gal(L/K)$-action extends to
$\Gamma$, this shows that the $\Gal(L/K)$-action extends to $\cal X$, as desired.
\qed\medskip

\begin{Remark}
  We stress that the reason for the extension of the $G$-action to ${\cal X}$ rather than ${\cal X}'$
  in the case of Abelian and hyperelliptic surfaces is their ``simple'' 
  geometry: they contain no smooth rational curves.
\end{Remark}

\subsection{The action on the special fiber}
In the situation of Proposition \ref{prop: extending group action to smooth model},
we now want  to understand whether the induced $G$-action on 
the special fiber ${\cal X}_0$ is trivial.
Quite generally, if $Y$ is a smooth and proper variety over some field $k$,
then the natural representation
$$
\begin{array}{ccccc}
  \rho_m &:& {\rm Aut}(Y) & \longrightarrow & {\rm Aut}\left(\,\Het{m}(Y_{\overline k},\QQ_\ell)\,\right) 
\end{array}
$$
is usually neither injective nor surjective.
We have the following exceptions:
\begin{enumerate}
  \item If $Y$ is an Abelian variety, then $\rho_1$ is injective.
    (Here, ${\rm Aut}(Y)$ denotes the automorphism group as an Abelian variety -- translations
     may act trivially on cohomology.)
  \item If $Y$ is a K3 surface, then $\rho_2$ is injective.
\end{enumerate}
Using these results (for references, see below), we have the following.

\begin{Proposition}
  \label{prop: trivial action on fiber}
  We keep the notations and assumptions of Proposition \ref{prop: rational group action}.
  Moreover, assume that either
  \begin{enumerate}
    \item $X$ is an Abelian surface and the $G$-action 
    on $\Het{1}(X_{\overline L},\QQ_\ell)$ is unramified, or 
    \item $X$ is a K3 surface and the $G$-action on $\Het{2}(X_{\overline L},\QQ_\ell)$ is unramified.
  \end{enumerate}
  Then, the $I_G$-action on $X_L$ extends to $\cal X$, and the induced $I_G$-action 
  on the special fiber ${\cal X}_0$ is trivial.
\end{Proposition}

\prf
We have already shown the extension of the $I_G$-action to ${\cal X}$ 
in Proposition \ref{prop: extending group action to smooth model}.
Moreover, the $I_G$-action on ${\cal X}_0$ is $k$-linear.
By assumption, the $I_G$-action on $\Het{m}(({\cal X}_0)_{\overline k},\QQ_\ell)$ is trivial
for $m=1,2$, respectively.

If $X$ is an Abelian surface, then the $I_G$-action on $({\cal X}_0)_{\overline{k}}$ is trivial
by the injectivity of $\rho_1$ in arbitrary characteristic, see, for example, 
\cite[Theorem 3 in Section 19]{Mumford Abelian}.
If $X$ is a K3 surface, then the $I_G$-action on $({\cal X}_0)_{\overline{k}}$ is trivial
by the injectivity of $\rho_2$ in arbitrary characteristic, see
\cite[Corollary 2.5]{Ogus} and \cite[Theorem 1.4]{Keum}
(in the case of complex and possibly non-algebraic K3 surfaces, 
see \cite[Proposition IX.6]{Asterisque} and \cite[Proposition 15.2.1]{Huybrechts}).
In both cases, the $I_G$-action on $({\cal X}_0)_{\overline{k}}$ is trivial,
and thus, also the original action on ${\cal X}_0$ is trivial.
\qed\medskip

\subsection{Tame quotients}
Now, in the situation of Proposition \ref{prop: extending group action to smooth model},
it is natural to study the quotient ${\cal X}/H$ and its special fiber, where $H$ is
a subgroup of $G$.
We start with the following easy result.

\begin{Proposition}
  \label{prop: tame quotient}
   Let $X$ be a smooth and proper variety over $K$.
   Let $L/K$ be a finite Galois extension with group $G$, such that
   $X_L$ admits a smooth model ${\cal X}\to\Spec\OO_L$.
   Moreover, assume that the natural $G$-action on $X_L$ extends
   to ${\cal X}$.
   Let $H$ be a subgroup of $G$ such that
   \begin{enumerate}
     \item $H$ is contained in the inertia subgroup of $G$,
     \item $H$ is of order prime to $p$, and
     \item $H$ acts trivially on the special fiber ${\cal X}_0$.
   \end{enumerate}
   Then,
    \begin{enumerate}
     \item the quotient ${\cal X}/H$ is smooth over $\Spec \OO_L^H$, and
     \item the special fiber of ${\cal X}/H$ is isomorphic to ${\cal X}_0$.
   \end{enumerate}
\end{Proposition}

\prf
First of all, the quotient ${\cal X}/H$ exists in the category of algebraic spaces
\cite[Chapter IV.1]{Knutson}.
Next, let $\widehat{{\cal X}}_0$ be the formal completion of ${\cal X}$ along the special fiber ${\cal X}_0$,
which is a formal scheme.
If $x\in{\cal X}_0$ is a closed point, then $\OO_{\widehat{{\cal X}}_0,x}$
is \'etale over the localization $A := \widehat{\OO}_L \langle y_1,y_2 \rangle_\idealm$
of the restricted power series ring 
$$
  \widehat{\OO}_L \langle y_1,y_2 \rangle \,:=\,
   \left\{  \left.
   \sum_{i_1, i_2 \geq 0}a_{i_1, i_2} y_1^{i_1} y_2^{i_2} \in \widehat{\OO}_L[[y_1,y_2]] \, \right\rvert \,
    \begin{array}{l} \ord_\pi(a_{i_1,i_2}) \to \infty \\ \mbox{ as } i_1+i_2 \to \infty \end{array}
 \right\}
$$
at the maximal ideal $\idealm = (\pi, y_1, y_2)$.
The induced $H$-action on the residue ring $A / (\pi) \iso \kappa(\OO_L)[y_1,y_2]_\idealm$ is trivial.
Thus, replacing $y_i$ by $\frac{1}{|H|}\sum_{\sigma\in H}\sigma(y_i)$ for $i=1,2$
(here, we use that the order of $H$ is prime to $p$) is simply a change of coordinates of $A$.
But then, the $H$-action on $A = \widehat{\OO}_L \langle y_1,y_2 \rangle_\idealm$ is trivial on $y_1$ and $y_2$,
and hence 
$\OO_{\widehat{({\cal X}/H)}_0,x} \iso \OO_{\widehat{{\cal X}}_0,x}^H$ 
is \'etale over $A^H \iso \widehat{\OO}_{L^H} \langle y_1,y_2 \rangle_\idealm$.
From this local and formal description, the smoothness of ${\cal X}/H$ follows immediately, and
we see that the quotient map ${\cal X}\to{\cal X}/H$ induces an isomorphism
of special fibers.
\qed\medskip

\subsection{Wild quotients}
Unfortunately, Proposition \ref{prop: tame quotient}
is no longer true if $p\neq0$ and 
$H$ is a subgroup of the inertia subgroup, whose order is divisible by $p$.
Let us illustrate this with a very instructive example.
We refer the interested reader to Wewers' article \cite{Wewers} for a more thorough treatment of wild 
actions and their quotients.

\begin{Example}
  Consider $K:=\QQ_p[\zeta_p]$, where 
  $\zeta_p$ is a primitive $p$.th root of unity.  
  Then, $\pi:=1-\zeta_p$ is a uniformizer in $\OO_K$, and the residue
  field is $\FF_p$, see \cite[Lemma 1.4]{Washington}, for example.
  Let $L$ be the finite extension $K[\varpi]$, where $\varpi:=\sqrt[p]{\pi}$.
  Then, $\varpi$ is a uniformizer in $\OO_L$, and the residue field is $\FF_p$, 
  that is, $L/K$ is totally ramified.
  By Kummer theory, $L/K$ is Galois  with
  group $H\iso\ZZ/p\ZZ$.
  More precisely, there exists a generator $\sigma\in H$ such that 
  $\sigma(\varpi)=\zeta_p\cdot \varpi$.
  We set 
  $$
     R\,:=\,\OO_L[x]
  $$
  and extend the $H$-action to $R$ by requiring
  that $\sigma(x)=\zeta_p^{p-1}\cdot x$.
  Then, we have $R/(\varpi)\iso\FF_p[x]$, and the induced $H$-action on the quotient
  is trivial.
  On the other hand, we find that
  $$
    R^H\,\iso\,\OO_K[x^p,x\cdot \varpi]\,\iso\,\OO_K[u,z]/(z^p-\pi u)
  $$
  is normal, but not regular  --
  this is an arithmetic version of the RDP singularity of type $A_{p-1}$.
  We also find that the special fiber
  $$
     R^H/ (\pi) \,\iso\,\FF_p[u,z]/(z^p)
  $$
  is not reduced.
  In particular, Proposition \ref{prop: tame quotient}
  does not extend to wild actions without extra assumptions.
  However, let us make two observations, whose significance will become clear
  in the proof of Proposition \ref{prop: wild quotient}.
  \begin{enumerate}
    \item The $H$-action on the special fiber $R/(\varpi)$ only {\em seems} 
        to be trivial, but in fact, it has become infinitesimal.
        More precisely, if $r\in R$ and $\overline{r}$ denotes its residue class
        in $R/(\varpi)$, then the $H$-action gives rise to a well-defined and non-trivial derivation 
  $$
  \begin{array}{cccccc}
    \theta &:& R/(\varpi) &\to& R/(\varpi) \\
    & & \overline{r} &\mapsto& \left(\, \frac{\sigma(r)-r}{\pi} \,\right) \mod \varpi.
  \end{array}
  $$
     \item The {\em augmentation ideal}, that is, the ideal of $R$ generated by all elements 
        of the form $\sigma(r)-r$ is not principal. 
        In fact, it can be generated by the two elements $\varpi^p x$ and $\varpi^{p+1}$.
    \end{enumerate}
\end{Example}
\medskip

Despite this example,
we have the following analog of Proposition \ref{prop: tame quotient}
in the wildly ramified case.
The main ideas of its proof are due to
Kir\'aly--L\"utkebohmert \cite[Theorem 2]{KL} and
Wewers \cite[Proposition 3.2]{Wewers}.

\begin{Proposition}
  \label{prop: wild quotient}
   Let $X$ be a smooth and proper variety over $K$.
   Let $L/K$ be a finite Galois extension with group $G$, such that
   $X_L$ admits a smooth model ${\cal X}\to\Spec\OO_L$.
   Moreover, assume that the natural $G$-action on $X_L$ extends
   to ${\cal X}$.
   Let $H$ be a subgroup of $G$ such that
   \begin{enumerate}
     \item $H$ is contained in the inertia subgroup of $G$,
     \item $H$ is cyclic of order $p\geq2$, and
     \item $H$ acts trivially on the special fiber ${\cal X}_0$.
   \end{enumerate}
   Then, the $H$-action induces a global and non-trivial derivation on 
   ${\cal X}_0$ or else both of the following two statements hold true
   \begin{enumerate}
       \item the quotient ${\cal X}/H$ is smooth over $\OO_L^H$, and
       \item the special fiber of ${\cal X}/H$ is isomorphic to ${\cal X}_0$.
   \end{enumerate}
\end{Proposition}

\prf
First of all, the quotient ${\cal X}/H$ exists in the category of algebraic spaces
\cite[Chapter IV.1]{Knutson}.
Next, we fix once and for all a generator $\sigma\in H$ and a uniformizer $\pi\in\OO_L$.
We use these to define the following: 
\begin{eqnarray*}
  N(\OO_L)      &:=& \max \left\{\, k\,\mid\, \pi^k \mbox{ divides } \sigma(x)-x \mbox{ for all } x\in\OO_L\,\right\} \\
  J_H (\OO_L) &:=& \mbox{ ideal of $\OO_L$ generated by $\sigma(x)-x$ for all $x\in \OO_L$. }
\end{eqnarray*}
Since $\OO_L$ is a DVR, the ideal $J_H(\OO_L)$ is principal.
More precisely, this ideal is generated by $y:=\sigma(\pi)-\pi$, and it is also generated
by $\pi^{N(\OO_L)}$.
In \cite{KL}, ideals generated by elements of the form $\sigma(x)-x$ are called {\em augmentation ideals}.
Also, it is not difficult to see that they do not depend on the choice of generator $\sigma$, which 
justifies the subscript $H$ rather than $\sigma$.

Next, let $\widehat{{\cal X}}_0$ be the formal completion of ${\cal X}$ along the special fiber ${\cal X}_0$,
which is a formal scheme.
For every point $x\in{\cal X}_0$, we define
\begin{eqnarray*}
  N(\OO_{\widehat{{\cal X}}_0,x})      &:=& \max \left\{\, k\,\mid\, \pi^k \mbox{ divides } \sigma(r)-r \mbox{ for all } r\in\OO_{\widehat{{\cal X}}_0,x}\,\right\} \\
  J_H (\OO_{\widehat{{\cal X}}_0,x}) &:=& \mbox{ ideal of $\OO_{\widehat{{\cal X}}_0,x}$ generated by $\sigma(r)-r$ for all $r\in \OO_{\widehat{{\cal X}}_0,x}$. }
\end{eqnarray*}
If $\eta\in{\cal X}_0$ denotes the generic point, then we have the following
\begin{equation}
 \label{eq: augmentation numbers}
   1\,\leq\,N(\OO_{\widehat{{\cal X}}_0,\eta}) \,\leq\, N(\OO_{\widehat{{\cal X}}_0,x}) \,\leq\, N(\OO_L),
\end{equation}
where the leftmost inequality follows from the triviality  of the $H$-action on ${\cal X}_0$.
We distinguish two cases:
\medskip

{Case (I):} $N(\OO_{\widehat{{\cal X}}_0,\eta})<N(\OO_L)$.

Let $x\in{\cal X}_0$ be an arbitrary point and set $R:=\OO_{\widehat{{\cal X}}_0,x}$
and $N_\eta:=N(\OO_{\widehat{{\cal X}}_0,\eta})$.
Then, we define a map
$$
\begin{array}{cccccc}
  \theta &:& R &\to& R/\pi R \\
  & & x &\mapsto& \left(\, \frac{\sigma(x)-x}{\pi^{N_\eta}} \,\right) \mod \pi,
\end{array}
$$
which is easily seen to be a derivation.
Since we have $N_\eta<N(\OO_L)$, we compute $\theta(\pi)=0$, and thus, 
$\theta$ induces a derivation
$\overline{\theta}:R/\pi R\to R/\pi R$.
This globalizes and gives rise to a derivation on the special fiber ${\cal X}_0$.
It follows from the definition of $N_\eta$ that this derivation is non-zero at
the generic point $\eta\in{\cal X}_0$, whence non-trivial.
\medskip

{Case (II):} $N(\OO_{\widehat{{\cal X}}_0,\eta})=N(\OO_L)$.

Let $x\in{\cal X}_0$ be an arbitrary point and set $R:=\OO_{\widehat{{\cal X}}_0,x}$.
Then, the two inequalities at the center and the right of (\ref{eq: augmentation numbers}) are equalities,
which implies that all inclusions in
$$
   \pi^{N(\OO_L)} \cdot R
   \,=\,J_H(\OO_L) \cdot R 
   \,\subseteq\, J_H(R) 
   \,\subseteq\, \pi^{N(R)}\cdot R 
$$
are equalities.
In particular,  $J_H(R)$ is a principal ideal, generated by $\pi^{N(\OO_L)}$.
But then, \cite[Proposition 5]{KL} implies
that there is an isomorphism
 of $R^H$- (resp. $\OO_L^H$-) modules
\begin{eqnarray*}
  R &\iso & R^H \,\oplus\,R^H\,\pi\,\oplus\,...\,\oplus\, R^H\,\pi^{p-1},\\
  \OO_L &\iso& \OO_L^H \,\oplus\,\OO_L^H\,\pi\,\oplus\,...\,\oplus\, \OO_L^H\,\pi^{p-1}.
\end{eqnarray*}
From this description, we conclude that the natural map
$$
  R^H \otimes_{\OO_L^H}\OO_L\,\to\,R
$$
is an isomorphism.
Moreover, if $\pi^H$ is a uniformizer of $\OO_L^H$, then the previous
isomorphism induces an isomorphism
$$
  R^H/\pi^H R^H \,\iso\, R/\pi R.
$$
This local computation at completions shows that ${\cal X}/H\times_{\OO_L^H}\OO_L$ 
is isomorphic to ${\cal X}$, and that the special fiber ${\cal X}_0$ 
of ${\cal X}$ is isomorphic to the special fiber of ${\cal X}/H$.
Since ${\cal X}$ is smooth over $\OO_L$, ${\cal X}_0$ is
smooth over the residue field of $\OO_L$, which implies that
also the special fiber of ${\cal X}/H$ is smooth over the residue
field of $\OO_L^H$.
But this implies that ${\cal X}/H$ is smooth over $\OO_L^H$.
\qed\medskip

Combining Propositions \ref{prop: tame quotient} and \ref{prop: wild quotient},
we obtain the following result.

\begin{Corollary}
  \label{cor: quotient}
   Let $X$ be a smooth and proper variety over $K$.
   Let $L/K$ be a finite Galois extension with group $G$, such that
   $X_L$ admits a smooth model ${\cal X}\to\Spec\OO_L$.
   Assume that the natural $G$-action on $X_L$ extends
   to ${\cal X}$, and that
   the inertia subgroup $I_G$ of $G$ acts trivially on the special fiber ${\cal X}_0$.
   Assume also that the special fiber ${\cal X}_0$ admits no non-trivial global vector fields.
   Then,
   \begin{enumerate}
       \item the quotient ${\cal X}/I_G$ is smooth over $\OO_L^{I_G}$, and
       \item the special fiber of ${\cal X}/I_G$ is isomorphic to ${\cal X}_0$.
   \end{enumerate}
\end{Corollary}

\prf
We have a short exact sequence
$$
   1\,\to\,P\,\to\,I_G\,\to\,T\,\to\,1,
$$
where $P$ is the unique $p$-Sylow subgroup of $I_G$, and
where $T$ is cyclic of order prime to $p$.
By definition, $P$ is the wild inertia, and $T$ is the tame inertia.

Being a $p$-group, $P$ can be written as a successive extension of cyclic groups of order $p$.
Thus, applying Proposition \ref{prop: wild quotient} inductively,  we obtain
a smooth algebraic space
$$
  {\cal X}/P \,\to\, \Spec\OO_{L}^P
$$
with special fiber ${\cal X}_0$, which is a model of $X_{L^P}$.

Finally, applying Proposition \ref{prop: tame quotient}
to the residual $T$-action on ${\cal X}/P$, we obtain a smooth algebraic space
$$
  {\cal X}/I_G \,\to\, \Spec\OO_L^{I_G}
$$
with special fiber ${\cal X}_0$, which is a model of $X_{L^{I_G}}$.
\qed\medskip

\begin{Remark}
\label{rem: RSK3}
If $X$ is a K3 surface, then the special fiber $\cal{X}_0$ is also a K3 surface and 
thus, admits no non-zero vector fields by a theorem of
Rudakov and Shafarevich \cite{rudakov shafarevich k3}.
\end{Remark}

\section{The N\'eron--Ogg--Shafarevich criterion}
\label{sec: main}

We now come to the main result of this article, which is
a criterion for good reduction of K3 surfaces, similar to
the classical N\'eron--Ogg--Shafarevich criterion for elliptic curves
and its generalization to Abelian varieties by Serre and Tate.
Then, we give a couple of corollaries concerning potential
good reduction, and good reduction after a tame extension.
Finally, we relate the reduction behavior of a polarized
K3 surface to that of its associated Kuga--Satake 
Abelian variety.

\subsection{The criterion}

Let us remind the reader of Section \ref{subsec: potential semi-stable reduction}, 
where we introduced Assumption $(\star)$ and established it in several 
cases.

\begin{Theorem}
  \label{thm: main}
  Let $X$ be a K3 surface over $K$ that satisfies $(\star)$.
  If the $G_K$-representation on $\Het{2}(X_{\overline{K}},\QQ_\ell)$ is
  unramified for some $\ell\neq p$, then
  \begin{enumerate}
    \item there exists a model of $X$ that is a projective scheme over $\OO_K$, 
       whose special fiber is a K3 surface with at worst RDP singularities.
    \item There exists an integer $N$, independent of $X$ and $K$, and a
       finite unramified extension $L/K$ of degree $\leq N$ such that
       $X_L$ has good reduction over $L$.
  \end{enumerate}
\end{Theorem}

\prf
By Theorem \ref{thm: matsumoto}, there exists a finite Galois extension $M/K$, 
say, with group $G$ and possibly ramified, such that there exists a smooth model of $X_M$
$$
 {\cal X}\,\to\,\Spec\OO_M \,.
$$
Choose an ample invertible sheaf ${\cal L}$ on $X$.
Then, by Proposition \ref{prop: termination of flops}, we can replace ${\cal X}$ by another smooth
model of $X$ such that the pull-back of ${\cal L}$ to $X_M$ restricts to an invertible sheaf on
${\cal X}_0$ that is big and nef.

Let $I_G$ be the inertia subgroup of $G$.
By Proposition \ref{prop: extending group action to smooth model}, the $I_G$-action
extends to ${\cal X}$ and by Proposition \ref{prop: trivial action on fiber},
the induced $I_G$-action on the special fiber ${\cal X}_0$ is trivial.
Thus, by Corollary \ref{cor: quotient} and Remark \ref{rem: RSK3},
the quotient 
$$
  {\cal X}/I_G \,\to\, \Spec\OO_L,
$$
where $L:=M^{I_G}$, is a model of $X_L$.
Since $L$ is a finite and unramified extension of $K$,
this establishes claim (2) except for the universal bound $N$.

The pull-back of ${\cal L}$ to ${\cal X}/I_G$ is still ample on the generic
fiber and big and nef when restricted to the special fiber.
By Proposition \ref{prop: rational group action},
there exists a birational morphism over $\Spec\OO_L$
$$
   \pi'\,:\, {\cal X}/I_G \,\to\, {\cal Y} ,
$$
that is an isomorphism on generic fibers, such that the geometric 
special fiber ${\cal Y}_0$ is a K3 surface with at worst RDP singularities, 
and such that the $H:={\rm Gal}(L/K)$-action on $X_L$ extends to ${\cal Y}$.
Since $L/K$ is unramified, the morphism
$\Spec\OO_L\to \Spec\OO_K$ is \'etale, from which it follows that the
quotient ${\cal Y}\to{\cal Y}/H$ is \'etale.
Thus, ${\cal Y}/H$ is a projective scheme over $\OO_K$, whose generic fiber is
$X$ and whose geometric special fiber is a K3 surface with at worst RDP 
singularities.
This establishes claim (1).

It remains to prove the existence of a universal bound $N$ in (2).
Since the Picard rank of a K3 surface is bounded above by $22$,
there is only a finite list $\sf L$ of Dynkin diagrams that is independent of the 
characteristic and whose associated root lattices 
can be embedded in the N\'eron--Severi lattice of a K3 surface.
Therefore, a K3 surface over an algebraically closed field has at most
$21$ RDP singularities, and all of them are from the list $\sf L$.
By \cite{Artin RDP}, there exist only finitely many analytic 
isomorphism types of RDP singularities with fixed dual resolution graph
over algebraically closed fields.
For every $k'$-rational singularity over some perfect field $k'$ that
becomes analytically isomorphic to a RDP singularity over $\overline{k}'$,
we have a versal deformation space ${\rm Def}$ 
over $k'$ or $W(k')$ (if ${\rm char}(k')=0$ or $>0$, respectively)
and a simultaneous resolution algebraic space ${\rm Res}$, 
which is finite over ${\rm Def}$ by \cite[Theorem 3]{Artin Brieskorn}.
Since deformation and resolution spaces solve universal problems,
the degree of ${\rm Res}\to{\rm Def}$ depends only on the analytic
isomorphism type of the singularity over $\overline{k}'$.
In particular, there exists an integer $N'$ such that every deformation of a
$k'$-rational singularity over $k'$ that becomes a RDP singularity
over $\overline{k}'$ from the list $\sf L$ can be resolved after an extension
of degree at most $N'$.
For each Dynkin diagram in $\sf L$ and for almost every characteristic
(in fact, it suffices to exclude $2 \leq p \leq 19$), 
there is only one analytic isomorphism type of RDP singularities and the 
corresponding degree of ${\rm Res}\to{\rm Def}$ is independent of the characteristic.
Therefore, the bound $N'$ can be taken to be independent of the characteristic.

Now, let ${\cal Y}_0$ be the special fiber of ${\cal Y}$.
Since $({\cal Y}_0)_{\overline{k}}$ has at most $21$ non-smooth points, 
all non-smooth points of ${\cal Y}_0$ become $k'$-rational after some finite 
extension $k'$ of $k$ of degree $\leq21!$.
Let $K'/K$ be the corresponding unramified extension of $K$.
From the previous discussion, it follows that after a (possibly ramified)
extension $L/K'$ of degree $\leq N'^{21}$, the surface 
$X_L$ has good reduction.
By the above arguments, we can descend a smooth model
of $X_L$ over $\OO_L$ to the the maximal unramified 
subextension $M$ of $L/K'$.
Thus, $X_M$ has good reduction, and $M/K$ is an unramified
extension of degree $\leq N:=21!\cdot N'^{21}$.
This establishes the bound claimed in (2).
\qed\medskip

\begin{Remark}
  In the statement (2) of Theorem \ref{thm: main}, we cannot avoid field extensions in general:
  in the next section, we will give examples of K3 surfaces $X$ over $\QQ_p$
  with unramified $G_{\QQ_p}$-representations
  on $\Het{2}(X_{\overline{\QQ}_p},\QQ_\ell)$ that do {\em not} admit
  smooth models over $\ZZ_p$.
\end{Remark}

\begin{Remark}
  Unlike curves and Abelian varieties, even if a K3 surface has good
  reduction over $L$, then a smooth model of $X_L$ over $\OO_L$
  need {\em not} be unique.
  However, by Proposition \ref{prop: flops along models},
  the special fibers of all smooth models are isomorphic and the models
  are connected by finite sequences of flopping contractions and their
  inverses (similar to the classical Atiyah flop).
\end{Remark}

If a smooth variety over $K$ has good reduction over an unramified extension,
then the $G_K$-representations on $\Het{m}(X_{\overline{K}},\QQ_\ell)$ are unramified
for all $m$ and for all $\ell\neq p$ by Theorem \ref{thm:obvious}. 
Thus, as in the case of Abelian varieties in \cite[Corollary 1 of Theorem 1]{serre tate}, we obtain
the following independence of the auxiliary prime $\ell$.

\begin{Corollary}
  \label{cor: independence of l}
  Let $X$ be a K3 surface over $K$ that satisfies $(\star)$. 
  Then, the $G_K$-representation on $\Het{2}(X_{\overline{K}},\QQ_\ell)$ is unramified
  for one $\ell\neq p$ if and only if it is unramified for all $\ell\neq p$.
\end{Corollary}

We remark that this independence of $\ell$ can be also derived from 
the weaker criterion of \cite{matsumoto2} (Theorem \ref{thm: matsumoto}),
combined with Ochiai's independence of traces \cite[Theorem B]{Ochiai}.

We leave the following easy consequence of Theorem \ref{thm: main}
to the reader
(this also can be derived from the criterion of \cite{matsumoto2}).

\begin{Corollary}
  Let $X$ be a K3 surface over $K$ such that the image of
  inertia
  $$
    \begin{array}{ccccc}
     \rho_\ell &:& I_K &\to& {\rm GL}\left(\, \Het{2}(X_{\overline{K}},\QQ_\ell)\,\right)
     \end{array}
  $$
  is finite. 
  If $X$ satisfies $(\star)$, then $X$ has potential good reduction.
 \qed
\end{Corollary}

If a $g$-dimensional Abelian variety over $K$
with $p>2g+1$ has potential good reduction,
then good reduction can be achieved over 
a {\em tame} extension of $K$ by \cite[Corollary 2 of Theorem 2]{serre tate}.
We have the following analog for K3 surfaces:

\begin{Corollary}
  \label{cor: tame extension}
  Let $X$ be a K3 surface over $K$ with potential good reduction.
  If $p\geq23$, then good reduction can be achieved after a
  tame extension.
\end{Corollary}

\prf
The idea of proof is the same as for Abelian varieties in \cite{serre tate}, we only
adjust the arguments to our situation:
since $X$ is projective, there exists an ample invertible sheaf $\cal L$ defined over $K$,
and then, its Chern class $c_1({\cal L})$ 
gives rise to a $G_K$-invariant class in $\Het{2}(X_{\overline{K}},\ZZ_\ell)(1)$.
Let $T^2_\ell$ be the orthogonal complement of $c_1({\cal L})$ 
with respect to the Poincar\'e duality pairing.
For $\ell\neq p$, we let
$$
   \rho_\ell  \,:\, G_{K} \,\to\, {\rm GL}\left(T^2_\ell\right)
$$
be the induced $\ell$-adic Galois representation, and denote by
$$
  {\rm red}_\ell \,:\, {\rm GL}\left( T^2_\ell \right) \,\to\,
  {\rm GL}\left( T^2_\ell / \ell  T^2_\ell \right) 
$$
its reduction modulo $\ell$.
As usual, we denote by $I_K$ (resp., $P_K$) the inertia (resp., wild inertia)
subgroup of $G_K$.
Since $X$ has potential good reduction, $\rho_\ell(I_K)$ is a finite group.
Moreover, if $\ell$ is odd, since $\ker{\rm red}_\ell$ has no non-trivial element of finite order
(as can be seen by taking the logarithm),
$\rho_\ell(I_K)$ is isomorphic to ${\rm red}_\ell\circ\rho_\ell(I_K)$ via ${\rm red}_\ell$.

Now, suppose that $\rho_\ell(P_K)$ is non-trivial.
Then, the order of ${\rm red}_\ell\circ\rho_\ell(I_K)$ is divisible by $p$ for all odd $\ell$.
In particular, if we set $n:={\rm rank}\,T_\ell^2=21$, then $p$ divides the order
$$
   \left| {\rm GL}_{n}(\FF_\ell) \right| \,=\, \ell^{n(n-1)/2}\,\prod_{s=1}^{n}(\ell^s-1)
$$
for all odd $\ell$.
By Dirichlet's theorem on arithmetic progressions, 
there exist infinitely many primes $\ell$ such that the residue class
of $\ell$ modulo $p$ generates the group $\FF_p^\times$, which is of order $p-1$.
Choosing such an $\ell$, we obtain the estimate $p-1 \leq n = 21$.
(When working directly with $H^2(X_{\overline{K}},\ZZ_\ell)(1)$ instead of the primitive cohomology group
$T^2_\ell$, we only get the estimate $p-1\leq22$, which includes the prime $p=23$.)
% note to ourselves: we also know that the image lies in the special orthogonal group, which yields the inequality 
% $p-1\leq20$ --- in particular, we do not gain anything: the proof is more complicated and we still only get $p\geq23$.

Thus, if $p\geq23$, then $\rho_\ell(P_K)$ is trivial.
But then, also the $P_K$-action on $\Het{2}(X_{\overline{K}},\ZZ_\ell)(1)$ is trivial.
Thus, there exists a {\em tame} extension $L/K$ such that the
$G_L$-action on $\Het{2}(X_{\overline{L}},\QQ_\ell)(1)$ is unramified.
By Theorem \ref{thm: main}, there exists an unramified extension of $L'/L$
such that $X_{L'}$ has good reduction.
In particular, $X$ has good reduction after a tame extension of $K$.
\qed\medskip

\subsection{Kuga--Satake varieties}
\label{subsec: Kuga Satake}

Given a polarized K3 surface $(X,{\cal L})$ over $\CC$,
Kuga and Satake \cite{Kuga Satake} associated to it
a polarized Abelian variety, the {\em Kuga--Satake Abelian variety}
$A:={\rm KS}(X,{\cal L})$, which is of dimension $2^{19}$.
Although their construction is transcendental, it is shown in work of
Rizov \cite{Rizov} and Madapusi Pera \cite{madapusi pera},
building on previous results of Deligne \cite{Deligne K3} and
Andr\'e \cite{Andre},
that the Kuga--Satake construction exists over arbitrary fields:
namely, if $(X,{\cal L})$ is a polarized K3 surface over some field $k$,
then ${\rm KS}(X,{\cal L})$ exists over some finite extension of $k$.
Then, we have the following relation between
good reduction of $(X,{\cal L})$ and
${\rm KS}(X,{\cal L})$.

\begin{Theorem}
  \label{thm: kuga-satake good reduction}
  Assume $p \neq 2$.
  Let $(X,{\cal L})$ be a polarized K3 surface over $K$.
   \begin{enumerate}
     \item If $X$ has good reduction, then ${\rm KS}(X,{\cal L})$
       can be defined over an unramified extension $L/K$, and it has good reduction over $L$.
     \item Assume that $X$ satisfies $(\star)$. Let
       $L/K$ be a field extension such that 
       both ${\rm KS}(X,{\cal L})$ and the Kuga--Satake correspondence (described below) can be 
       defined over $L$.
       If ${\rm KS}(X,{\cal L})$ has good reduction over $L$, then $X$ has good reduction
       over an unramified extension of $L$.
  \end{enumerate}
\end{Theorem}

\prf
We will use the notations and definitions of \cite{madapusi pera}.

(1)
The pair $(X,{\cal L})$ gives rise to a morphism
$\Spec K\to{\mathsf M}_{2d}^\circ$, where ${\mathsf M}_{2d}^\circ$ denotes
the moduli space of primitively polarized K3 surfaces of degree $2d:={\cal L}^2$.
By assumption, there exists a smooth model of $X$ over $\OO_K$, and 
by Proposition \ref{prop: termination of flops}, there even exists a smooth
model ${\cal X}$ of $X$ over $\OO_K$, such that the restriction of $\cal L$
to the special fiber is big and nef.
Thus, the morphism $\Spec K\to{\mathsf M}_{2d}^\circ$ extends to a morphism
$\Spec\OO_K\to{\mathsf M}_{2d}$, where ${\mathsf M}_{2d}$ denotes the moduli space
of primitively quasi-polarized K3 surfaces of degree $2d$.
Passing to an unramified extension $L/K$ of degree $\leq2$ if necessary,
the previous classifying morphism extends to a morphism
$\Spec\OO_L\to\widetilde{{\mathsf M}}_{2d}$, see \cite[Section 5]{madapusi pera}.
Composing with the morphism 
$\widetilde{{\mathsf M}}_{2d}\to{\cal S}(\Lambda_d)$ from \cite[Proposition 5.7]{madapusi pera},
we obtain a morphism $\Spec\OO_L\to{\cal S}(\Lambda_d)$.
We recall from \cite[Section 4]{madapusi pera} that there exists a finite and \'etale cover
$\widetilde{{\cal S}}(\Lambda_d)\to{\cal S}(\Lambda_d)$, such that the Kuga--Satake
Abelian scheme is a relative Abelian scheme over $\widetilde{{\cal S}}(\Lambda_d)$.
Thus, after replacing $L$ by a finite and unramified extension if necessary, we can lift
the latter morphism to a morphism $\Spec\OO_L\to\widetilde{{\cal S}}(\Lambda_d)$.
Thus, we obtain a Kuga--Satake Abelian variety ${\rm KS}(X,{\cal L})$ over $L$
that has good reduction, and where $L$ is an unramified extension of $K$.

(2)
By assumption, ${\rm KS}(X,{\cal L})$ is defined over $L$ and has good reduction over $L$.
Thus, the $G_L$-action on $\Het{1}({\rm KS}(X,{\cal L})_{\overline{L}},\QQ_\ell)$ is unramified.
By the usual properties of the Kuga--Satake construction, there exists a
%$G_L$-equivariant embedding
embedding
$$
   \Pet{2}(X_{\overline{K}},\QQ_\ell)(1) \to {\rm End}\left(
    \Het{1}({\rm KS}(X,{\cal L})_{\overline{L}},\QQ_\ell)
   \right),
$$
where $\Pet{2}(X_{\overline{K}},\QQ_\ell)(1)$ denotes the orthogonal complement of 
$c_1({\cal L})$ inside $\Het{2}(X_{\overline{K}},\QQ_\ell)(1)$,
and this embedding is $G_L$-equivariant by assumption.
This implies that also the $G_L$-action on $\Pet{2}(X_{\overline{K}},\QQ_\ell)(1)$
is unramified.
Since ${\cal L}$ is defined over $K$, the $G_L$-action on the $\QQ_\ell$-subvector
space generated by $c_1({\cal L})$ inside $\Het{2}(X_{\overline{K}},\QQ_\ell)(1)$ 
is trivial.
From this, we conclude that the $G_L$-action on $\Het{2}(X_{\overline{K}},\QQ_\ell)(1)$
is unramified.
By Theorem \ref{thm: main}, $X$ has good reduction over an unramified extension of $L$.
\qed\medskip

\begin{Remark}
  Let us make two comments:
  \begin{enumerate}
  \item If $(X,{\cal L})$ is a polarized K3 surface with good reduction, then the previous
  theorem asserts that ${\rm KS}(X,{\cal L})$ can be defined over an unramified 
  extension $L$ of $K$. 
  Thus, if ${\rm KS}(X,{\cal L})$ can be descended to some field $K'$ 
  with $K\subseteq K'\subseteq L$ (so far, not much is known about fields of definition
  of Kuga--Satake Abelian varieties), then,
  since $L/K'$ is unramified and since ${\rm KS}(X,{\cal L})$ has good reduction over $L$
  by the previous theorem, the descended Abelian variety will have good
  reduction over $K'$ by \cite{serre tate}.
  \item We can almost remove the $p\neq2$ 
   hypothesis in Theorem \ref{thm: kuga-satake good reduction}:
   by \cite[Proposition A.12]{kim madapusi pera} 
   (see also the proof of \cite[Theorem A.1]{kim madapusi pera}), 
   there exists a Kuga--Satake morphism with the properties needed to make
   the proof of Theorem \ref{thm: kuga-satake good reduction} work
   also in residue characteristic 2, but so far only outside the locus 
   of superspecial K3 surfaces.
  \end{enumerate}
\end{Remark}

\section{counterexamples}
\label{sec: counterexample}

In this final section we give examples of K3 surfaces $X$ over  $\QQ_p$ for all 
$p \geq 5$ with unramified $G_{\QQ_p}$-representation on 
$\Het{2}(X_{\overline{\QQ}_p},\QQ_\ell)$ that do {\em not} have good reduction over $\QQ_p$.
In particular, the unramified extension from Theorem \ref{thm: main} needed
to obtain good reduction may be non-trivial.
The examples in question already appeared in \cite[Section 5.3]{matsumoto2}
and rest on examples due to van~Luijk \cite[Section 3]{van luijk}.

\begin{Example}
   \label{example: no good reduction}
  Fix a prime $p\geq5$.
  We choose integers $a,c$ such that
  $a\not\equiv 0,\frac{27}{16}\mod p$, such that $c\equiv1\mod8$,
  and such that $c$ is not a quadratic residue modulo $p$.
  Then, we choose a homogeneous polynomial $f\in\ZZ[x,y,z,w]$
  of degree $3$, such that the following congruences hold
  $$
      \begin{array}{ccll}
         f &\equiv& \phi & \mod 2\\
         f &\equiv& x^3+y^3+z^3+aw^3 &\mod p,
      \end{array}
  $$
  where $\phi$ is as in Example \ref{example: resolutions}.
  Finally, we define the quartic hypersurface
   $$
    {\cal X}\,:=\,{\cal X}(p) \,:=\, \left\{ wf+ \left(pz^2+xy+\frac{p}{2}yz\right)^2-\frac{cp^2}{4}y^2z^2\,=\,0
    \right\} \,\subset\, \PP^3_{\ZZ_p},
  $$
  and denote by $X=X(p)$ its generic fiber.
 \end{Example}
 
 \begin{Theorem}
   Let $p\geq5$ and let $X$ and ${\cal X}$  be
   as in Example \ref{example: no good reduction}.
   Then, $X$ is a smooth K3 surface over $\QQ_p$, such that
   \begin{enumerate} 
      \item the $G_{\QQ_p}$-representation on $\Het{2}(X_{\overline{\QQ}_p},\QQ_\ell)$ is unramified for all $\ell\neq p$,
      \item ${\cal X}$ is a projective model of $X$ over $\ZZ_p$, 
         whose geometric special fiber is a K3 surface with RDP singularities of type $A_1$,
       \item $X$ has good reduction over the unramified extension
         $\QQ_p[\sqrt{c}]$,
       \item $X$ does not have good reduction over $\QQ_p$.
   \end{enumerate}
 \end{Theorem}
 
 \prf
 Smoothness of $X$ follows from considering the equations over $\ZZ$, 
 reducing modulo $2$ and checking smoothness there. 
 Claims (2) and (3) are straightforward computations (for (3), blow up the ideal ${\cal I}_+$ or ${\cal I}_-$ defined below), and since $X$ has
 good reduction after an unramified extension, also claim (1) follows.
 We refer to \cite[Section 5.3]{matsumoto2} for computations and details.
 
To show claim (4), we argue by contradiction, and assume that there exists
a smooth and proper algebraic space ${\cal Z}\to\Spec\ZZ_p$ with generic
fiber $X$.
Since the generic fibers of ${\cal X}$ and ${\cal Z}$ are isomorphic, 
such an isomorphism extends to a birational, but possibly rational map
$$
   \alpha\,:\,{\cal Z}\,\dashrightarrow\,{\cal X}.
$$
Next, let ${\cal L}$ be an ample invertible sheaf on ${\cal X}$,
for example, the restriction of $\OO(1)$ from the ambient ${\PP^3_{\ZZ_p}}$.
Restricting $\cal L$ to the generic fiber ${\cal X}_\eta$,  and pulling it back via $\alpha$, 
we obtain an ample invertible sheaf $\alpha_\eta^*({\cal L}_\eta)$ on ${\cal Z}_\eta$.
Since ${\cal Z}$ is smooth over $\Spec\OO_K$, this invertible sheaf on ${\cal Z}_\eta$
extends uniquely to an invertible sheaf on ${\cal Z}$ that we denote by ${\cal M}$.

By Proposition \ref{prop: termination of flops}, there exists a rational and birational map
$$
   \varphi\,:\,{\cal Z} \,\dashrightarrow\, {\cal Z}^+,
$$
where ${\cal Z}^+$ is another model of $X$ with good reduction, and such that
the transform ${\cal M}^+$ of $\cal M$ on ${\cal Z}^+$ is ample on the generic
fiber, and big and nef on the special fiber.
We denote by $\alpha^+:{\cal Z}^+\dashrightarrow{\cal X}$ the composition
$\alpha\circ\varphi^{-1}$.
Then, 
$$
 {\cal Z}^+\,\to\,({\cal Z}^+)'\,:=\,\Proj\bigoplus_{n\geq0} H^0\left({\cal Z}^+,({\cal M}^+)^{\otimes n}\right)
$$
is a birational morphism that contracts precisely those curves on the special fiber
${\cal Z}^+_0$ that have zero-intersection with ${\cal M}_0^+$, and nothing else.
By construction, we have $(\alpha^+_\eta)^*{\cal L}\iso{\cal M}^+$ and thus,
by \cite[Theorem 5.14]{Kovacs}, there exists an isomorphism
$({\cal Z}^+)' \stackrel{\sim}{\to} {\cal X}$ over $\Spec\ZZ_p$.

Thus, we have shown that the model ${\cal X}$ admits a simultaneous
resolution $\alpha^+ \colon {\cal Z}^+ \to {\cal X}$ of singularities over $\ZZ_p$.
But then, let $x \in {\cal X}_0$ be an $\FF_p$-rational singular point,
for example, the point $x = w = y+z = 0$.
Then, let $\OO_{{\cal X},\overline{x}}$ be the strict local ring, and denote by
${\rm Cl}(\OO_{{\cal X},\overline{x}})$ its Picard group.
Then, $\alpha^+$ induces a $G_{\QQ_p}$-equivariant surjection
$(R^1 \alpha^+_* \OO_{{\cal Z}^+}^*)_{\bar x} \to {\rm Cl}(\OO_{{\cal X},\overline{x}})$. 
However, this is impossible for the following reason:
\begin{enumerate}
 \item The group $(R^1 \alpha^+_* \OO_{{\cal Z}^+}^*)_{\bar x}$ is generated by the class of the exceptional curve, 
    which is $\FF_p$-rational, and thus the $G_{\QQ_p}$-action on it is trivial.
 \item The $G_{\QQ_p}$-action on ${\rm Cl}(\OO_{{\cal X},\overline{x}})$ is non-trivial.
    More precisely, if we define the following ideals of $\OO_{{\cal X},\overline{x}}$
    $$
       {\cal I}_{\pm}:=\left(w, \left(pz^2+xy+\frac{p}{2}yz\right) \pm \frac{\sqrt{c}}{2}p yz \right),
    $$
    then their classes in ${\rm Cl}(\OO_{{\cal X},\overline{x}})$ satisfy $[{\cal I}_+] = - [{\cal I}_-] \neq [{\cal I}_-]$, and since
    $G_{\QQ_p}$ acts on $\QQ_p[\sqrt{c}]$ as $\sqrt{c}\mapsto -\sqrt{c}$,
    the $G_{\QQ_p}$-action on ${\rm Cl}(\OO_{{\cal X},\overline{x}})$ is non-trivial.
 \end{enumerate}   
This contradiction shows that $X$ does not have good reduction over $\QQ_p$,
and establishes claim (4).
\qed\medskip

\begin{Remark} 
 \label{rem:phi}
  This example is the one the second named author gave in \cite[Section 5.3]{matsumoto2}.
  There, the choice of $\phi$ ensured that $X$ was a smooth K3 surface, as well of 
  Picard number one.
  In the present paper, we only need smoothness, and therefore,
  could have used a simpler polynomial.
  The same remark holds for Example \ref{example: resolutions}.
\end{Remark}

\end{document}